\documentclass[a4paper]{amsart}
\usepackage{amsmath,amsthm,amssymb,mathtools, stmaryrd,
	tikz-cd,xcolor}
\usepackage[cmtip]{xypic}
\usepackage[T1]{fontenc}

\usepackage[colorlinks=true,allcolors=blue!60!black]{hyperref}
\usepackage[capitalise]{cleveref}
\usepackage{xifthen}

\usepackage[left=3.5cm,right=3.5cm,top=3.25cm,bottom=3.25cm]{geometry}
\usepackage{setspace}
\theoremstyle{plain}
\newtheorem{thm}{Theorem}[section]

\newtheorem{prop}[thm]{Proposition} 
\newtheorem{lem}[thm]{Lemma}
\newtheorem{cor}[thm]{Corollary}
\newtheorem{con}{Condition}
\theoremstyle{definition} 
\newtheorem{defn}[thm]{Definition} 
\newtheorem{exmp}[thm]{Example}
\theoremstyle{remark}
\newtheorem{rem}[thm]{Remark}

\newcommand{\Z}{\mathbb{Z}}

\DeclareMathOperator{\MF}{MF}
\DeclareMathOperator{\Cof}{Cof}
\DeclareMathOperator{\CofII}{Cof^{\mathrm{II}}}
\DeclareMathOperator{\D}{D}
\DeclareMathOperator{\DII}{D^{{II}}}
\DeclareMathOperator{\DIIdg}{D^{II}_{{dg}}}
\DeclareMathOperator{\Dco}{D^{co}}
\DeclareMathOperator{\Dctr}{D^{ctr}}
\DeclareMathOperator{\Dctrdg}{D^{ctr}_{dg}}
\DeclareMathOperator{\DCoh}{DCoh}
\DeclareMathOperator{\End}{End}
\DeclareMathOperator{\HH}{HH}

\DeclareMathOperator{\HHc}{\mathcal{HC}}

\DeclareMathOperator{\HHctr}{\mathcal{HC}_{\mathrm{ctr}}}
\DeclareMathOperator{\HHII}{HH^{\mathrm{II}}}
\DeclareMathOperator{\HHIIc}{\mathcal{HC}^{\mathrm{II}}}
\DeclareMathOperator{\Ho}{Ho}
\DeclareMathOperator{\Hom}{Hom}
\DeclareMathOperator{\Homdg}{\underline{Hom}}
\DeclareMathOperator{\MC}{MC}
\DeclareMathOperator{\Perf}{Perf}
\DeclareMathOperator{\PerfII}{Perf^{\mathrm{II}}}
\DeclareMathOperator{\LC}{LC}
\DeclareMathOperator{\RHom}{RHom}
\newcommand{\RHomII}{\mathrm{RHom}^{\mathrm{II}}}
\DeclareMathOperator{\tw}{tw}
\DeclareMathOperator{\Tw}{Tw}
\newcommand*\lon{\nobreak \mskip6mu plus1mu \mathpunct{} \nonscript \mkern-\thinmuskip{:}\mskip2mu \relax}
\newcommand{\op}{\mathrm{op}}
\newcommand{\cat}[1]{\mathsf{#1}}
\newcommand{\cop}{\mathsf{op}}
\newcommand{\DGVect}{\cat{DGVect}}
\newcommand{\PreVect}{\cat{PreVect}}
\newcommand{\cuAlg}{\cat{cuAlg}}
\newcommand{\id}{\operatorname{id}}

\newcommand{\dgm}{\cat{Mod}\textrm{-}}
\newcommand{\DGMod}[2][]{%
	\ifthenelse{\equal{#1}{}}{\dgm #2}{#1\textrm{-}\dgm #2}%
} 
\newcommand{\DGModI}[2][]{%
	\ifthenelse{
		\equal{#1}{}}{(\dgm #2)^{\mathrm{I}}}{(#1\textrm{-}\dgm #2)^{\mathrm{I}}
	}%
} 
\newcommand{\DGModII}[2][]{%
	\ifthenelse{\equal{#1}{}}{(\dgm #2)^{\mathrm{II}}}{(#1\textrm{-}\dgm #2)^{\mathrm{II}}}%
} 
\newcommand{\DGModIIi}[2][]{%
	\ifthenelse{\equal{#1}{}}{(\dgm #2)^{\mathrm{II}}_{\mathrm{c, inj}}}{(#1\textrm{-}\dgm #2)^{\mathrm{II}}_{\mathrm{c, inj}}}%
} 

\DeclareFontFamily{U}{mathx}{\hyphenchar\font45} 
\DeclareFontShape{U}{mathx}{m}{n}{<-> mathx10}{}
\DeclareSymbolFont{mathx}{U}{mathx}{m}{n}
\DeclareMathAccent{\widecheck}{0}{mathx}{"71} 
\DeclareMathAccent{\widebar}{0}{mathx}{"73} 

\newcommand{\cB}{\widecheck{\mathsf{B}}}
\newcommand{\hB}{\widehat{\mathsf{B}}}
\numberwithin{equation}{section}

\reversemarginpar
\newif\ifshowcomments
\showcommentstrue 

\theoremstyle{plain}
\newtheorem{innercustomthm}{Theorem} 
\newenvironment{customthm}[1]
{\renewcommand\theinnercustomthm{#1}\innercustomthm}
{\endinnercustomthm}

\newtheorem{innercustomcor}{Corollary} 
\newenvironment{customcor}[1]
{\renewcommand\theinnercustomcor{#1}\innercustomcor}
{\endinnercustomcor}

\newtheorem{innercustomdef}{Definition} 
\newenvironment{customdef}[1]
{\renewcommand\theinnercustomdef{#1}\innercustomdef}
{\endinnercustomdef}
\title[Hochschild cohomology of the second kind]{Hochschild cohomology of the second kind:\\ Koszul duality and Morita invariance}
\author{A. ~Guan}
\address{School of Mathematics\\
	University of Birmingham\\
	Birmingham B15 2TT\\United Kingdom}
\email{a.guan@bham.ac.uk}

\author{J. Holstein}
\address{Department of Mathematics\\
	Universit\"at Hamburg\\
	20146 Hamburg\\
	Germany	
}
\email{julian.holstein@uni-hamburg.de}

\author{A.~Lazarev}

\address{Department of Mathematics and Statistics\\
	Lancaster University\\
	Lancaster LA1 4YF\\United Kingdom}
\email{a.lazarev@lancaster.ac.uk}

\author{K. Rasmussen}
\address{Department of Mathematics\\
	Universit\"at Hamburg\\
	20146 Hamburg\\
	Germany	
}
\email{kristoffer.rasmussen@uni-hamburg.de}

\thanks{This work was partially supported by EPSRC grant EP/T029455/1. 
	The second and fourth author acknowledge support by the Deutsche Forschungsgemeinschaft (DFG, German Research Foundation) through EXC 2121 ``Quantum Universe'' -- project number 390833306.  The second author also acknowledgs support through SFB 1624 ``Higher structures, moduli spaces and integrability'' -- project number 506632645.  }

\begin{document}
	
	\bibliographystyle{../hsiam2}
	\begin{abstract}
		We define Hochschild cohomology of the  second kind for differential graded (dg) or curved algebras as a derived functor in the twisted derived category.
		The Hochschild cohomology of the second kind of a dg or curved algebra $A$	is then equivalent to the classical Hochschild cohomology of the twisted derived dg category of $A$, which is often geometrically meaningful.
		Examples include the category of $\infty$-local systems on a topological space, the bounded derived category of a complex manifold and the category of matrix factorizations. 
		
		We also show that Hochschild cohomology of the second kind is preserved under (nonconilpotent) Koszul duality and weak equivalences of curved algebras.
		The main technical ingredient is a new bimodule version of Koszul duality.
	\end{abstract}
	
	\maketitle
	
	\tableofcontents
	
	\section{Introduction}
	Hochschild cohomology of an associative algebra (or, more generally, a dg algebra or dg category) is an important invariant that has been of much current interest. 
	It is known to be invariant under Morita equivalence: given any (dg) algebra its Hochschild cohomology is isomorphic to the Hochschild cohomology of its category of perfect dg modules. 
	A related fundamental result, or series of results, is the invariance of Hochschild cohomology under Koszul duality. In other words, Hochschild cohomology of a dg algebra and that of its Koszul dual (dg or curved) conilpotent coalgebra are isomorphic, and the isomorphism preserves various structures that they possess; we refer to \cite{kel21a} for the most structured version of this result and its history. 
	
	There are several situations when classical Hochschild cohomology, also called the Hochschild cohomology of the first kind, does not capture the information of interest.
	This is the case if one works with a curved algebra, or if one is interested in studying exotic versions of the derived category, different derived categories of the second kind.
	
	One answer to this is compactly supported Hochschild cohomology, also called Hochschild cohomology of the second kind \cite{pp12, caldararu2013curved}.
	Note that the  Hochschild complex of a dg algebra is a \emph{double complex} and for the Hochschild cohomology of the first kind one takes its direct products totalization whereas for the Hochschild cohomology of the second kind the totalization is by the direct sums. For an ungraded algebra that makes no difference, but for a genuine dg algebra this leads to completely different theories, the one of the second kind being much more mysterious and less understood. 
	Hochschild cohomology of the second kind belongs to `homological algebra of the second kind', where the notion of a quasi-isomorphism is replaced by a (typically finer) notion of a weak equivalence. 
	For the general philosophy of homological algebra of the second kind we refer to Positselski's survey \cite{Positselski23}. 

	One key motivation to study these new invariants is that derived categories of the second kind carry important geometric information in many examples.
	For example, the singular cochain algebra on a topological space, the de Rham algebra of a smooth manifold or the Dolbeault algebra of a complex analytic manifold have twisted derived categories (a certain type of derived categories of the second kind) that capture much more geometric information than their ordinary derived categories, i.e.\ the category of $\infty$-local systems and the bounded derived category of coherent sheaves respectively \cite{chl21, blo10}. 
	One now expects that Hochschild cohomology of the second kind for these algebras is also geometrically meaningful and in particular computes the ordinary Hochschild cohomology of the associated dg categories. This indeed happens in some situations e.g. for categories of matrix factorizations, \cite{pp12}. In contrast, the ordinary Hochschild cohomology of these algebras are much less interesting.

	 In the present paper, we will construct for any curved algebra $A$ a new Hochschild cohomology of the second kind $\HHIIc(A)$ which recovers the ordinary Hochschild cohomology of the twisted derived category.
	 Our construction, which arises naturally as a derived functor of the second kind, does not, in general, agree with compactly supported Hochschild cohomology; at the same time, it has a number of advantages over the latter: it is Morita invariant (in a suitable sense) and compatible with Koszul duality. We give an explicit complex representing $\HHIIc(A)$ and compute $\HHIIc(A)$  in geometrically interesting examples related to $\infty$-local systems, coherent sheaves and matrix factorizations.

	In more detail, our results are as follows: To set the scene, we recall in Section  \ref{sect:background} the notion of a twisted complex over a curved category (in particular over a curved algebra or dg category).
	
	Then we recall the twisted derived category $\DII(A)$ for a dg (or curved) algebra $A$ which is compactly generated by its subcategory $\PerfII(A)$ (which just the closure of finitely generated twisted complexes under homotopy idempotents).
	The category $\DII(A)$ has also been called the compactly generated derived category of the second kind, see Remark \ref{R:name}.
	
	We recall 
	that category $\DII(A)$ is equivalent to the contraderived category $\Dctr(\cB A)$ of its Koszul dual, where $\cB A$ is the \emph{extended} bar-construction of $A$, cf.\ \cite{aj13, gl21, bl23} concerning this notion.
	
		Note that Hochschild cohomology, be it of the first or second kind, is constructed as a (co/contra)derived functor of bi(co)modules; whereas Koszul duality is usually formulated as an equivalence between one-sided modules and comodules. It is natural, therefore, to establish Koszul duality as an equivalence between bimodules over an algebra and bicomodules over a suitably Koszul dual coalgebra;  this is also the approach of \cite{kel21a}.
Our generalization to	Koszul duality for contraderived and twisted derived categories of bimodules is established in Sections \ref{sect:cobar} and \ref{sect:bimkoszul}.
%
	
	 There is a significant subtlety here since for a curved algebra $A$ one has to distinguish between the twisted derived category of $A$-bimodules, $\DII(\DGMod[A]{A})$ and the twisted derived category of $A^\op \otimes A$-modules.
	They have the same underlying ordinary category which carries two different natural model structures. 
	In Theorem \ref{T:twbim} we establish the model structure for $A$-bimodules and then go on to show it is naturally Koszul dual to category of bimodules over the pseudocompact algebra $\cB(A^\op)\otimes \cB(A)$.
This is the category in which Hochschild cohomology arises as a derived endomorphism algebra. 

	Our key definition is then the following:
	\begin{customdef}{\ref{D:main}}\label{def:HochII}
		The \emph{Hochschild cohomology complex of the  second kind} of a dg or curved algebra $A$ is
		$\HHIIc(A) = \RHom_{\DII(\DGMod[A]{A})}(A,A)$ and the \emph{Hochschild cohomology of the second kind} is $\HHII(A) = H(\HHIIc(A))$. 
	\end{customdef}


	We construct a complex computing this invariant in Corollary \ref{C:compute} by replacing the usual bar construction by an extended bar construction in the standard Hocshchild cochain complex.
	We then have our main theorem:

	\begin{customthm}{\ref{thm:hhperf2}}
		For any curved algebra $A$ there are quasi-isomorphisms of dg algebras	$\HHIIc(A) \simeq \HHc(\PerfII(A)) \simeq \HHc(\DII(A))$.
	\end{customthm}

	This has a number of consequences, in particular that Hochschild cohomology of the second kind is invariant under MC equivalences and naturally carries a Gerstenhaber algebra structure.
	
	In Section \ref{sect:koszulinvariance} we show that Hochschild cohomology of the second kind is equivalent to Hochschild cohomology of its Koszul dual pseudocompact algebra:
	\begin{customcor}{\ref{C:computeB}}\label{cor:HHII}
		For any dg algebra $A$, there is a quasi-isomorphisms of dg algebras
		\[
		\HHIIc(A) \simeq \HHctr(\cB A).
		\]
	\end{customcor}
	From this, a Morita invariance result for coHochschild cohomology of coalgebras also follows.
	
	We consider Hochschild cohomology of the second kind for geometric and topological examples in Section \ref{sect:examples}:
		\begin{itemize}
		\item If $A$ is the Dolbeault algebra of a smooth complex manifold $X$ then $\HHII(A)$ is isomorphic the Hochschild cohomology of the bounded derived category of coherent sheaves on $X$. If $X$ is algebraic, this is computed by cohomology of the exterior algebra on vector fields on $X$.
		\item If $A$ is the singular cochain algebra of a topological space $M$ or the de Rham algebra of a smooth compact manifold $M$ then $\HHII(A)$ is isomorphic to the Hochschild cohomology of the category of $\infty$-local systems on $M$,
		\item If $R_w$ is $\Z/2$-graded regular commutative algebra concentrated in even degrees with curvature element $w$ then $\HHII(R_w)$ is isomorphic to the Hochschild cohomology of the category of matrix factorizations $\MF(R,w)$ which encodes the structure of the hypersurface singularity defined by $w$. It can be computed by an explicit complex of derivations twisted by $dw$.
	\end{itemize}

	As was already mentioned, our version of Hochschild cohomology of the second kind is \emph{not} the same as that of Positselski and Polishchuk \cite{pp12} but rather lies between the ordinary Hochschild cohomology and that of Positselski-Polishchuk. 
	In fact, our Hochschild cochain complex of the second kind is a certain completion of the one considered by Positselski-Polishchuk and the ordinary Hochschild cochain complex is a further completion, cf.\ Remark \ref{remark:comparison} for a precise statement.

	It is worth emphasizing that for our version of Hochschild cohomology of the second kind, the statement of Theorem \ref{thm:hhperf2}  holds with no assumptions on the curved algebra $A$. 
	In contrast, the typical assumption needed for this kind of Morita invariance for Positselski-Polishchuk's version is a kind of ``resolution of the diagonal'' property. We consider the closely related condition  that the natural functor 
	$$
	i \colon \Perf(\PerfII(A^\op)\otimes \PerfII(E))\to  \PerfII(A^\op\otimes E)
	$$
	that sends a pair $M, N$ of perfect modules of the second kind to $M \otimes N$
	is an equivalence.
	This is not true in general; instead, we only have that 
	$	\PerfII(A^\op \otimes E) \simeq \PerfII(\PerfII(A^\op) \otimes \PerfII(E)),$ see Proposition \ref{L:perfperf}. 
	
	If the condition above is satisfied, we may identify the twisted derived category of $(A,E)$-bimodules with the twisted derived category of $A^\op\otimes E$-modules.
	If, moreover, certain further conditions on the underlying graded algebra hold, then the Hochschild cohomology constructed by Polishchuk-Positselski agrees with our construction, see Theorem \ref{thm:comp HH_PP and HH2}.
	
	In particular, this applies to the example of matrix factorizations (assuming 0 is the only critical value of $w$), giving a conceptual explanation why the two constructions (with very different underlying complexes) define the same Hochschild cohomology of second kind.

	\subsection{Acknowledgements} The authors would like to thank Patrick Antweiler for pointing out  mistakes in previous versions of this paper. ChatGPT 5.6 was used for some classical algebra.
	
	We thank the referee for several useful comments.
	\section{Background}\label{sect:background}
	\subsection{Notation and conventions}
	We fix a perfect ground field $k$ throughout. We also fix a grading group which we assume to be $\Z$ for purposes of exposition, but which may be taken to be $\Z/2$. If the grading group is $\Z/2$ we also assume $\operatorname{char}(k) \neq 2$.
	
	 The category of (cohomologically graded) differential graded (dg) vector spaces over $k$ will be denoted by $\DGVect$. 
	 The category of graded vector spaces with a derivation (that need not square to zero) is denoted by $\PreVect$.
	 
	 The shift of a dg vector space $A$ is denoted by $A[1]$ so that $A[1]^i=A^{i+1}$. The category $\DGVect$ is monoidal with respect to the tensor product; monoids and comonoids in it are dg algebras and dg coalgebras respectively. 	
	 A \emph{curved} algebra $(A, d, h)$ is a graded algebra $A$ equipped with an element $h \in A^2$ and a derivation $d$ (sometimes called a \emph{pre-differenial}) such that $d^2(a) = [h,a]$. If $h=0$, a curved algebra can be viewed as a dg algebra. 
	 Note however, that a morphism of curved algebras from $(A, d, h)$ to $(A', d', h')$ takes the from $(f,a)$ where $f: A \to A'$ is a map of graded algebra
	 and $a \in A'$ is an element of degree one such that $f(dx)= d'f(x)+[a, f(x)]$ and $f(h) = h' + d'(b)+b^2$. Thus, dg algebras do not form a full subcategory of curved algebras.
	 
	 A left curved module $(M,d)$ over $(A,d,h)$ is a graded left $A$-module $M$ with a derivation $d$ such that $d^2(m) = hm$. For right curved modules, the latter identity is modified as $d^2(m) = -mh$.
	 We will work with left and right modules, but the default will be right modules.
	 
	 Given (curved) algebras $A$ and $E$, we will often consider $(A,E)$-bimodules which  will be denoted as $\DGMod[A]{E} := \DGMod{(A^\op \otimes E)}$. 
	 
	 The notions of a curved coalgebra and a comodule over it can be defined by duality.
	 More details on curved (co)algebras and (co)modules can be found e.g.\ in \cite{pos11}. 
	 
	 The underlying graded of a dg or curved algebra or module $A$ is denoted by $A^\#$.
	
	\subsection{Pseudocompact algebras and modules} 
	A pseudocompact (pc) graded vector space is a projective limit of finite-dimensional graded vector spaces. A natural example of a pc graded vector space is $V^*:=\Hom(V,k)$, the $k$-linear dual to a discrete graded vector space $V$. As such, it is equipped with the topology of the projective limit and the continuous linear dual $\Hom(V^*,k)$ gives us $V$ back. Thus, the category of pc graded vector spaces is anti-equivalent to the category of discrete graded vector spaces under the functor of $k$-linear duality. Similarly, a pc graded curved algebra is a projective limit of finite-dimensional graded curved algebras; a pc graded curved module over a pc graded curved algebra is a projective limit of finite dimensional graded curved modules. A pc curved graded algebra is always the dual to a graded curved coalgebra; this establishes an anti-equivalence between pc graded curved algebras and graded curved coalgebras. Note that conilpotent (curved) coalgebras correspond to pc augmented local (curved) algebras.  Similarly the category of pc graded curved modules over a pc graded curved algebra is anti-equivalent to the category of graded curved comodules over the corresponding graded curved coalgebra.

	A tensor product of pc algebras or modules is always assumed completed,  equivalently it is the linear dual of the tensor product of the dual coalgebras or comodules.

	\subsection{Curved categories} 
	A \emph{curved category} $\cat{A}$ is a category enriched over $\PreVect$  such that each object has a curvature element.
	To be precise, for each pair of objects $X,Y \in \cat{A}$, homomorphisms from $X$ to $Y$ form a $k$-precomplex which we denote by $\Homdg_\cat{A}(X,Y)$ with pre-differential $d$. For each pair $f:X\to Y$ and $g:Y\to Z$ of composable morphisms, the equation $d(gf) = d(g) f+(-1)^{|g|}g d(f)$ holds.
	Furthermore, for each object $X\in \cat{A}$, there is a 2-cocycle $h_X\in \Homdg_\cat{A}(X,X)$ satisfying $d^2(f)=h_Yf-fh_X$ for all $f\in \Homdg_\cat{A}(X,Y)$. If $h_X=0$ for all $X \in \cat{A}$, we say that $\cat{A}$ is a \emph{dg category}.
	
	A \emph{curved functor} is a pair $(F,a)$, with $F:\cat{A}\to \cat{B}$ being an additive functor and $a$ a collection of fixed degree 1 elements $a_X\in 
	\Homdg_\cat{B}(F(X),F(X))$, such that $F(d(f))=d(F(f))+a_YF(f)-(-1)^{|f|}F(f)a_X$ and $F(h_X)=h_{F(X)}+d(a_X)+a^2_X$.
	A curved functor is said to be \emph{strict} if $a_X=0$ for all $X\in \cat{A}$. A \emph{dg functor} is always a strict curved functor between dg categories.
	
	Any dg category $\cat{A}$ has an associated ordinary category $\Ho\cat{A}$, called its \emph{homotopy category}, whose objects are the same as $\cat{A}$ but whose morphisms are defined by $\Ho\cat{A}(X,Y) = H^0(\Homdg_{\cat{A}}(X,Y))$.
	
	A dg functor $F \colon \cat{A} \to \cat{B}$ is a \emph{quasi-equivalence} if
	\begin{enumerate}
		\item $F$ induces quasi-isomorphisms $\cat{A}(X,Y) \to \cat{B}(FX,FY)$ for any objects $X, Y$ in $\cat{A}$. 
		\item $H^0(F)$ is essentially surjective (or, equivalently assuming (1), $H^0(F)$ is an equivalence of categories). 
	\end{enumerate}
	In this paper, we are concerned with dg categories and curved algebras, however for ease of exposition it will be useful to have an ambient category containing both notions.
	\subsection{Bar and cobar constructions} Let $A$ be a dg algebra; we choose a section of the unit map $k\to A$ so that $A\cong k\oplus \bar{A}$; if $A$ is augmented then this can be chosen as a splitting of dg algebras but in general it need not be. 
	Then the \emph{bar-construction} $\operatorname{B}(A)$ of $A$ is the cofree conilpotent coalgebra $T\bar{A}[1]$; the differential and multiplication on $A$ determine the structure of a curved coalgebra on $T\bar{A}[1]$.    Similarly, for a curved coalgebra $C$ with a choice of splitting of a counit map so that $C\cong k\oplus \bar C$, its \emph{cobar-construction} $\Omega(C)$ is the free algebra $T\bar{C}[-1]$; the differential and comultiplication on $C$ determine the structure of a curved algebra on $\Omega(C)$. We will refer to \cite[Section 6.1]{pos11} for detailed definitions and properties of the bar and cobar constructions.
	
	In this paper we dualize all coalgebras and work with equivalent data of pc algebras. According to this principle, we consider the dual bar-construction $\hB(A):=\hat{T}\bar{A^*}[-1]$ where $\hat{T}$ stands for the completed tensor algebra; this is a curved local pc algebra (that, of course, carries the same information as the curved coalgebra $\operatorname{B}(A)$). 
	
	Similarly, instead of a curved coalgebra $C$, we will work with the curved pc algebra $D:=C^*$; then by definition, its cobar-construction $\Omega(D)$ is the curved algebra $\Omega(D):=T\bar{D}^*[-1]$, isomorphic to $\Omega(C)$. Note that the functors $\Omega$ and $\hB$ thus defined are contravariant.
	
	Finally, we need the notion of an \emph{extended bar-construction}. Given a dg or curved algebra $A$, the underlying coalgebra of the extended bar-construction of $A$, is the cofree \emph{nonconilpotent} coalgebra on $\bar{A}[1]$; according to our convention we dualize it and work instead with the pc curved algebra $\cB A$; the pc completion of the free algebra $T\bar{A}^*[-1]$. The differential and curvature are induced by muliplication, differential and curvature of $A$. 
	Details of this construction are given in \cite[Definition 2.5]{gl21}.

	
	\subsection{Curved modules and twistings} 
	Let $\cat{A}$ be a curved category. A (\emph{right}) \emph{curved $\cat{A}$-module} is a strict curved functor $M \colon \cat{A}^{\op} \to \PreVect$.
	
	We denote by $\DGMod \cat{A}$ the dg category of curved $\cat{A}$-modules. 
	Similarly a \emph{left} curved $\cat A$-module is a right curved $\cat A^\op$-module, i.e.\ a strict curved functor $N: \cat A \to \PreVect$. We will always specify if we consider a left dg module; by default a dg/curved module will be a right module.
	
	If $\cat{A}$ is a dg category, the notion of a curved $\cat{A}$-module coincide with the usual notion of a dg $\cat{A}$-module, i.e. a dg functor $M\colon \cat{A}^\cop\to \DGVect$.
	A dg $\cat{A}$-module $M$ is said to be \emph{acyclic} if it is acyclic pointwise, i.e.\ $M(X)$ is an acyclic complex for all $X \in \cat{A}$. 
	The \emph{derived category} $\D(\cat{A})$ of dg $\cat{A}$-modules is the Verdier quotient of $\Ho (\DGMod \cat{A})$ by acyclic dg $\cat{A}$-modules.
	There is a model structure, called the projective model structure, on $\DGMod \cat{A}$, where weak equivalences are pointwise quasi-isomorphisms and fibrations are pointwise surjections.
	All objects are fibrant; we denote by $\Cof(\cat A)$ the cofibrant objects in $\DGMod \cat{A}$.
	By general results on model categories, $\Ho \Cof(\cat{A}) = \D(\cat{A})$.
	
	Twisted complexes for dg categories were first defined in \cite{bk91} and later redefined in \cite{bll04}.
	
	\begin{defn}
		A (\emph{two-sided}) \emph{twisted complex} over $\cat{A}$ is a formal expression $(\bigoplus_{i=1}^n C_i[r_i], q)$, where $C_i \in \cat{A}$, $r_i \in \mathbb{Z}$, $n \geq 0$, $q=(q_{ij})$, $q_{ij} \in \Hom(C_j[r_j], C_i[r_i])$ homogeneous of degree 1 such that $(dq+q^2)_{ij}$ is $h_{C_i}$ if $i=j$ and 0 otherwise. 
		A twisted complex is \emph{one-sided} if $q_{ij} = 0$ for all $i \ge j$.
		
		For two twisted complexes $C$ and $C'$, the space of \emph{morphisms of twisted complexes} $\Hom(C,C')$ is the $\mathbb{Z}$-graded $k$-module of matrices $f=(f_{ij})$, $f_{ij} \in \Hom(C_j[r_j], C'_i[r'_i])$ with differential 
		$df = (df_{ij}) + q'f -(-1)^{|f_{ij}|}fq$. 
		Composition of morphisms is usual matrix multiplication.
	\end{defn}

	We denote the dg category of twisted complexes over $\cat{A}$ by $\Tw(\cat{A})$ and the dg category of one-sided twisted complexes by $\tw(\cat{A})$.
	In \cite{bk91}, these were respectively denoted $\text{Pre-Tr}(\cat{A})$ and $\text{Pre-Tr}^+\!(\cat{A})$; furthermore, $\tw(\cat A) = \text{Pre-Tr}^+\!(\cat{A})$ can alternatively be defined as the pretriangulated hull of $\cat{A}$, that is, it is the closure of $\cat{A}$ under shifts and cones. 
	
	There is a natural functor from twisted complexes over a dg category $\cat{A}$ to right dg $\cat{A}$-modules which sends $(\bigoplus_{i=1}^n C_i[r_i], q)$ to $(\bigoplus_{i=1}^n \Hom_{\cat{A}}(-, C_i)[r_i], q_*)$.
	
	As an example, consider the case where $\cat{A}$ is a curved algebra $(A,d,h)$, considered as a dg category with one object.
	Recall that a \emph{Maurer--Cartan element} in $A$ is an element $x \in A^1$ such that $h + dx+x^2 = 0$, and that the set of all Maurer--Cartan elements in $A$ is denoted by $\MC(A)$.
	Then a twisted complex over $A$ is a pair $(M, q)$ where $M \cong V \otimes A$ as an $A$-module for some finite-dimensional $\mathbb{Z}$-graded vector space $V$, and $q \in \MC(\End V \otimes A)$.
	The $A$-module $M$ becomes a (right) dg $A$-module when equipped with the differential $1 \otimes d_A + q$, and in fact, every curved $A$-module structure on the $A^\#$-module $M$ arises this way, as noted in \cite[Remark 3.2]{chl21}.

	A twisted complex over $A$ is therefore precisely the same data as a \emph{finitely generated} twisted $A$-module in the following sense.

	\begin{defn}\label{twisted module}
		A \emph{twisted module} over a curved algebra $(A, d, h)$ is a curved $A$-module whose underlying graded module has the form $V \otimes A^\#$ for  a graded $k$-module $V$.
		Explicitly, a twisted module is a pre-differential graded module $(V \otimes A, 1 \otimes d + q)$ equipped with the natural $A$-action where  re-writing $q$ as an element in $(\End V \otimes A)^1$ it satisfies $dq+q^2 = 1 \otimes h \in \End V \otimes A$.	
	\end{defn}
If  $A$ is curved then $A$ itself is not a (left or right) twisted module over itself; it is, however, always a bimodule over itself.

	There is a slight clash of terminology here; nevertheless this is unavoidable as we wish to refer to non-finitely generated twisted modules later.
	Note also that $(V\otimes A, 1 \otimes d_A)$ in Definition \ref{twisted module} is a curved $(\End V \otimes A)$-$A$-\emph{bi}module whose differential has been \emph{twisted} by the element $q \in \MC(\End V \otimes A)$. 
	More generally, we have the following.
	
	\begin{defn}\label{def:twisted}
		Let $(A,d,h)$ be a curved algebra and $x \in A^1$.
		\begin{enumerate}
			\item The \emph{twisted algebra of $A$ by $x$}, denoted $A^x = (A,d^x,h^x)$, is the curved algebra with the same underlying algebra as $A$, differential $d^x(a) = d(a) + [x,a]$ and curvature element $h^x=h+d(x)+x^2$. If $x\in \MC(A)$, the curvature element of $A^x$ is $0$.
			\item Let $(M,d_M)$ be a \emph{left} dg $A$-module.
			The \emph{twisted module of $M$ by $x$}, denoted $M^{[x]} = (M,d^{[x]})$, is the left dg $A^x$-module with the same underlying module structure as $M$ and differential $d^{[x]} (m) = d(m) + xm$.
		\end{enumerate}
		 If $M$ is a curved $A$-$B$-bimodule for curved algebras $A$ and $B$ and $x\in A^1$ then $M^{[x]}$ is a curved $A^x$-$B$-bimodule. 
	\end{defn}
	
	A \emph{perfect} dg $\cat{A}$-module for a dg category $\cat{A}$ is a cofibrant dg $\cat{A}$-module that is homotopy equivalent to a direct summand of a module in $\tw(\cat{A})$.
	We denote by $\Perf(\cat{A})$ the full dg subcategory of $\Cof(\cat{A})$ consisting of perfect modules.
	It can be shown that $\D(\cat{A})$ is a \emph{compactly generated} triangulated category, and that $\Ho \Perf(\cat{A})$ consists precisely of the compact objects in $\D(\cat{A})$, cf.\ \cite{Toen07a}. 
	
	\subsection{DG Morita equivalence}
	
	A dg functor $F \colon \cat{A} \to \cat{B}$ is a (\emph{dg}) \emph{Morita equivalence} if the induced map $F_! \colon \D(\cat{A}) \to \D(\cat{B})$ is an equivalence of triangulated categories.
	
	A main result of \cite{bk91} is that the Yoneda embedding $h \colon \cat{A} \to \tw(\cat{A}) \hookrightarrow \Tw(\cat{A})$ induces dg equivalences $\tw(\cat{A}) \to \tw(\tw(\cat{A}))$ and $\Tw(\cat{A}) \to \Tw(\Tw(\cat{A}))$.
	In particular, this implies that $h_! \colon \Perf(\cat{A}) \to \Perf(\Perf(\cat{A}))$ is a quasi-equivalence, so the Yoneda embedding is a Morita equivalence.
	In fact, $\Perf(\cat{A})$ is a Morita fibrant replacement of $\cat{A}$ in the sense of \cite[Section~4.6]{kel06a}.
	
\section{Koszul duality}\label{sect:koszul}
	
	\subsection{Model structures of the  second kind and Koszul duality}
	
	In this section we introduce some model structures on (co)modules that will feature in the remainder of this paper.
	These model structures will all give rise to ``derived categories of the  second kind'', the collective name given to derived categories which arise from localizing at some collection of weak equivalences that are finer than quasi-isomorphism.
	
	There are many such derived categories; we will recall the contraderived category of a pc curved algebra, following \cite{pos11}, and the twisted derived category of a curved algebra, introduced under the name ``compactly generated derived category of the second kind'' in \cite{gl21}.
	These are compatible via Koszul duality, which we will also recall.
	
	In the next subsection we will introdue the twisted derived category of bimodules.
	
	We first define, following \cite{pos11}, that in a category of dg or curved modules or comodules an object is \emph{coacyclic}
	(resp.\ \emph{contraacyclic}) if it lies in the smallest subcategory containing all totalizations of exact triples and closed under cones and direct sums (resp.\ direct products).

	Let now $A$ be curved algebra and $C$ be a pc curved algebra over a field $k$. 
	The categories $\DGMod{C}$ and $\DGMod{A}$, of pc curved $C$-modules and curved $A$-modules respectively, admit the following model category structures of the second kind.
	
	\begin{thm}[see {\cite[Theorem 8.2]{pos11}}] \label{T:cmcPC}
		For any pc curved algebra $C$, the category $\DGMod{C}$ is a model category where a morphism $f \colon M \to N$ is
		\begin{enumerate}
			\item a \emph{weak equivalence} if the cone of the dual map $f^* \colon N^* \to M^*$ of $C^*$-comodules is coacyclic; 
			\item a \emph{fibration} if it is surjective;
			\item a \emph{cofibration} if it has the left lifting property with respect to acyclic fibrations.
		\end{enumerate}
	\end{thm}
	Equivalently to the dual having a coacyclic cone, we can also characterize weak equivalences directly as having a contraacyclic cone.
	
\begin{defn}	For a pc curved algebra $C$, the homotopy category of $\DGMod{C}$ is the \emph{contraderived category} of $C$, denoted by $\Dctr(C)$.
	\end{defn}
		
	\begin{thm}[see {\cite[Theorem 4.6]{gl21}}] \label{T:cmcII}
		Let\/ $A$ be a curved algebra.
		There is a model category structure on\/ $\DGMod A$, where a morphism $f \colon M \to N$ is
		\begin{enumerate}
			\item a \emph{weak equivalence} if it induces a quasi-isomorphism
			\[
			\Homdg_A(T, M) \to \Homdg_A(T, N)
			\]
			for any finitely generated twisted $A$-module $T$;
			\item a \emph{fibration} if it is surjective;
			\item a \emph{cofibration} if it has the left lifting property with respect to acyclic fibrations.
		\end{enumerate}
	\end{thm}
	
	Where it is necessary to distinguish between the usual model structure of the first kind on dg $A$-modules (where weak equivalences are quasi-isomorphisms) and the model structure of \cref{T:cmcII}, we will denote these model categories respectively by $\DGModI{A}$ and $\DGModII{A}$.
	We say that a dg $A$-module is \emph{cofibrant of the  second kind} if it is cofibrant in $\DGModII{A}$, and denote these cofibrant objects by $\CofII(A)$.

\begin{defn}
	We denote the homotopy category of $\DGModII A$ by $\DII(A)$ and call it the \emph{twisted derived category} of $A$.
	\end{defn}	
	\begin{rem}\label{R:name}
		In previous work (including earlier versions of this paper) the category $\DII(A)$ was called the  \emph{(compactly generated) derived category of the  second kind}, which we now consider too unwieldy. 
		We thus propose the name twisted derived category as it is based on twisted complexes, and is indeed equivalent to the derived category of the category of twisted complexes over $A$, see Theorem \ref{T:compactII} below. 
		Note that despite the name this is not obtained from an (ordinary) derived category by a twisting procedure.
		
		The name twisted derived category is also sometimes used for the derived category of twisted sheaves on an algebraic variety, see \cite{caldararu2000derived}.
		It turns out that (at least on a smooth complex variety) the twisted derived category in the sense of Caldararu is an example of a twisted derived category in our sense, namely the twisted derived category of a curved version of  the Dolbeault algebra, see \cite[Setion 7.4]{antweiler}.				
	\end{rem}
	The model structures above are then related by the following statement, referred to as Koszul duality for the cobar construction.
	
	\begin{thm}\label{T:KdCobar}
		Let $C$ be a pc curved algebra. There is a Quillen equivalence
		\[
		G \colon (\DGMod{C})^\cop \rightleftarrows \DGModII{\Omega C} \lon F.
		\] 
	\end{thm}

		We will sometimes write $G_C(N)$ for $GN$ when we want to make the dependence on $C$ clear.
		
	The functor $G$ sends a pc curved $C$-module $N$ to the curved $\Omega C$-module 
	\[
	GN := (N^*\otimes \Omega C)^{[\xi]},
	\]
where $\xi$ is the canonical Maurer--Cartan element in $C \otimes \Omega C$. The twist makes sense as we regard $N^* \otimes \Omega C$ as a $(C \otimes \Omega C)^{\op} \otimes \Omega C$-module,
	 where the left action of $C$ comes from the right action of $C$ on $N$ and the $\Omega C^\op \otimes \Omega C$-action arises from the left and right multiplication on $\Omega C$. 
	Hence the twisted module $(N^*\otimes \Omega C)^{[\xi]}$ is a left $(C \otimes \Omega C)^{\xi}$-module and a right $\Omega C$-module.
	As the twisted algebra $(C \otimes \Omega C)^{\xi}$ has 0 curvature we may forget the left $(C \otimes \Omega C)^{\xi}$-action to obtain the right $C$-module
	$GN$.
	More explicitly, $GN$ can be written as a precomplex as follows after choosing a $k$-linear splitting $C\cong \widebar{C}\oplus k$:
	\[
	N^* \leftarrow N^* \otimes \widebar{C}^* \leftarrow N^* \otimes (\widebar{C}^*)^2 \leftarrow \dots.
	\]
	This is, in fact, a double complex where the vertical pre-differential (not indicated explicitly) 
	is induced by the internal pre-differentials in $C$ and $N$ whereas the horizontal pre-differential is induced by the multiplication in $C$ and its action on $N^*$. 
	Explicitly, the horizontal pre-differential is given by 
	\[
	d(c_0\otimes\ldots \otimes c_k\otimes n)
	=\sum_{i=0}^{k-1}(-1)^ic_0\otimes\ldots c_{i}c_{i+1}\ldots c_k\otimes n +(-1)^{k+1}c_0\otimes\ldots \otimes c_kn
	\]
	where the last term in the above sum corresponds to the twisting by $\xi$. 
	Similarly, the functor $F$ sends a curved $\Omega C$-module $M$ to the pc curved $C$-module
	\[
	FM := (M^* \otimes C)^{[\xi]}
	\]
	where we use the left $\Omega C$-module structure on $M^*$ and the left $C$-module structure on $C$ to define the twist.
	
	\begin{rem}
		A version of Theorem~\ref{T:KdCobar}, formulated in the language of coalgebras and comodules, is due to \cite[Theorem 6.7]{pos11}. Combining Theorem \ref{T:KdCobar} with Positselski's result we see that if $A= \Omega C$ then the homotopy category $\DII(A)$ is precisely the coderived category of $A$ defined by localizing at coacyclic $A$-modules.
		
		In the case that $C$ is augmented and local (equivalently, that $C^*$ is conilpotent), the cobar construction $\Omega C$ is cofibrant and one recovers the more classical statement where $\DGModII{\Omega C}$ on the right hand side replaved with the category of $\Omega C$-modules with the ordinary projective model structure.
		
		More on the relationship between different (co)derived categories can be found in \cite[Section 3.3]{gl21}.
	\end{rem}
	
	For the sake of convenience, we will now sketch the proof of \cref{T:KdCobar}. 
	First note that the functors $F$ and $G$ are adjoint in the sense that there are natural isomorphisms
	\[
	\Hom_{\DGMod\Omega C}(GN, M) \cong \Hom_{\DGMod C}(FM,N),
	\]
	for $N$ and $M$ as above, as both sides are identified with the dg vector spaces $(N\otimes M)^{[\xi]}$. 
	Composing the functors gives the pc $C$-module
	\[
	FGN = (N \otimes \Omega C^* \otimes C)^{[\xi\otimes 1+1\otimes \xi]},  
	\]
	which as a (double) pre-complex, can be written as follows:
	\[
	N \otimes C \leftarrow N \otimes \widebar{C} \otimes C 
	\leftarrow N \otimes (\widebar{C})^2 \otimes C \leftarrow \dots.
	\]
	The horizontal pre-differential is described by a similar formula as above.
	Note, however, that due to the tensor factor $C$, the horizontal pre-differential is acyclic except at the first term; in fact it gives the standard cobar resolution of $N$, the map from the above complex to $N$ being induced by the action map $N\otimes C\to N$. 
	Of course, one has to be careful calling this a resolution since $C$ is a curved pc algebra and $N$ is a pc $C$-module.
	This can be made sense of as follows. 
	Consider the homotopy cofiber of the map $FGN \to N$; this is the double complex that can be represented as follows:
	\[
	N \leftarrow N \otimes C \leftarrow N \otimes \widebar{C} \otimes C 
	\leftarrow N \otimes \widebar{C}^2 \otimes C \leftarrow \dots.
	\]
	Its horizontal pre-differential is in fact a differential and moreover acyclic (forgetting about differentials and curvatures this is just the bar resolution of $N$).
	By taking canonical truncations this complex can be represented
	as a homotopy inverse limit of acyclic complexes of finite horizontal length, which are, therefore, contra-acyclic as pc $C$-modules. 
	Thus, we obtain the following result, from which \cref{T:KdCobar} follows.
	
	\begin{prop}\label{P:resolutiononesided}
		The adjunction unit $FGN \to N$ is a cofibrant resolution of the pc dg $C$-module $N$.
	\end{prop}

	Analogous statements can be obtained if we instead start with a curved algebra $A$.
	
	If $A$ is a dg algebra, depending on whether we consider model structures of the first or second kind over $\DGMod{A}$, we will have two different equivalences. Note that there is the ordinary (dual) bar-construction $\hB A$ of $A$ and the extended (dual) bar-construction $\cB A$, both viewed as pc curved algebras; 
	however, the latter may be used for a curved algebra $A$.

	We then have the following Koszul duality statements for the bar constructions, see \cite[Theorem 4.7]{gl21}:
	
	\begin{thm}\label{T:KdBar}
		Let $A$ be a dg algebra there are Quillen equivalences
		\[
		(\DGMod{\hB A})^\cop \rightleftarrows \DGModI{A}
		\]
		and
		\[
		(\DGMod{\cB A})^\cop \rightleftarrows \DGModII{A}.
		\]
		and the second equivalence also holds for $A$ a curved algebra.
	\end{thm}
	
	We collect some useful facts about the model structure of the  second kind:
	\begin{prop}\label{prop:useful}
		\begin{enumerate}
			\item Twisted modules are cofibrant of the second kind if  they can be represented as a filtered colimit along cofibrations of finite generated twisted modules. All cofibrant modules are retracts of such.
			\item\label{it:coacyclic} Any map $f: M \to N$ in $\DGMod A$ with coacyclic or contra-acyclic cone is also a weak equivalence in the sense of Theorem \ref{T:cmcPC}.
		\end{enumerate}
	\end{prop}
	\begin{proof}
		\begin{enumerate}
			\item Any $\cB A$-module is fibrant, and may be written as a filtered limit of finite-dimensional modules along fibrations. Under Koszul duality the finite-dimensional $\cB A$-modules correspond to finitely generated twisted $A$-modules. An arbitrary cofibrant object is a retract of its cofibrant replacement (which is obtained as a union of finitely generated twisted modules).
			\item
			Let $0 \to L \to M \to N \to 0$ be an exact triple and $T$ a finitely generated twisted $A$-module. Then $0 \to \Hom_A(T, L) \to \Hom_A(T, M) \to \Hom_A(T, N) \to 0$ is exact. Therefore the dg space of homomorphisms from $T$ into the totalization of $L \to M \to N$ is acyclic. Thus the weak equivalences are closed under totalization of exact triples. They are automatically closed under direct products, and as any finitely generated twisted module is compact they are also closed under direct sums.\qedhere
		\end{enumerate}
	\end{proof}
	
	We now define perfect modules of the  second kind. 
	
	\begin{defn}
		A curved $A$-module is \emph{perfect of the  second kind} if it is homotopy equivalent to a direct summand of a two-sided twisted module in $\Tw(A)$.
		We denote by $\PerfII(A)$ the full dg subcategory of $\CofII(A)$ consisting of perfect modules of the  second kind. 
	\end{defn}
	We note that any homomorphism $f: A\to B$ induces a dg functor $f_!=\otimes_A B: \PerfII(A)\to\PerfII(B)$, where $B$ is endowed with a left $A$-module structure induced by $f$. 

	\begin{thm} \label{T:compactII}
		Let $A$ be a curved algebra. 
		Then $\Ho \PerfII(A)$ is the category of compact objects in $\DII(A)$, and $\DII(A)$ is compactly generated. Furthermore, there is an equivalence of triangulated categories $\DII(A) \cong D(\PerfII(A))$.
	\end{thm}
	
	\begin{proof}
		By definition, $\Ho \PerfII(A)$ consists of curved $A$-modules that are direct summands of finitely generated twisted $A$-modules, and these are compact (as image of finite pseudocompact modules under Koszul duality).
		
		To prove that $\PerfII(A)$ generates, we use that $\DII(A)$ is by Koszul duality (Theorem \ref{T:KdBar}) anti-equivalent as a triangulated category to the contraderived category of pc  dg modules over a pc  dg algebra, which is cocompactly cogenerated by finite dimensional dg modules by \cite[Section 5.5]{pos11} (which correspond to finitely generated twisted $A$-modules).
		
	To show that $\PerfII(A)$ contains all compact objects we use \cite[Lemma 2.2]{neeman1992connection} in the case $\mathcal S = \mathcal R = \DII(A)$.
	
	For the last part we consider the natural Yoneda embedding $\DGMod
 A \to \DGMod{\PerfII(A)}$ restricted to fibrant-cofibrant objects in $\DGMod A$. 
 This sends homotopy equivalences to homotopy equivalences and thus induces a functor on triangulated categories. 
 It identifies the subcategories $\Tw(A)$, which are compact generators, and preserves coproducts (as elements in $\Tw(A)$ are compact in the dg category $\DGMod A$). Thus, it is an equivalence \cite[Lemma 3.3]{schwede2001stable}. 
	\end{proof}
	\subsection{Model structure of curved algebras}
	
	The twisted derived category is generally not a quasi-isomorphism invariant of dg algebras and for curved algebras the notion of quasi-isomorphism does not make sense as they have no cohomology. However, there exists a model structure on the category $\cuAlg$ of curved algebras whose weak equivalences are described as follows:
	
	\begin{defn}\label{D:mc}
		Let $A$ and $B$ be curved algebras.
		A map $f \colon A \to B$ is a \emph{Maurer-Cartan equivalence} (MC equivalence for short) if for any curved coalgebra $C$ the induced map between MC moduli sets $\mathcal{MC}(C,A) \to \mathcal{MC}(C,A)$ is a bijection.
		Here $\mathcal{MC}(C,A)$ is the set of isomorphism classes of objects in the homotopy category of the dg category of MC elements in the convolution algebra $\Hom(C,A)$.
	\end{defn}
	The twisted derived category is an MC equivalence invariant.
	\begin{prop}[\cite{bl23}]\label{prop:mcd2}
		If $f:A\to B$ is an MC equivalence of curved algebras, then
		$$f_!:\DII(A)\to \DII(B)$$
		is a triangle equivalence
	\end{prop}
	
	\begin{exmp}\label{exmp:cobarbar}
		For any curved algebra $A$ the counit map $\Omega \cB A \to A$ is an MC equivalence. This is \cite[Theorem 6.4]{bl23}. 
	\end{exmp}

	We write $\cuAlg$ for the model category with MC equivalences. There is a "dual" model structure on the category of curved coalgebras making $(\Omega,\cB)$ a Quillen pair. See \cite[Section 6]{bl23} for a more detailed discussion of these model structures.
	\subsection{Twisted derived categories of bimodules}
	
	In homological algebra of the first kind to consider a derived category of $(A, E)$-bimodules it suffices to consider the derived category of $A^\op\otimes E$-modules.
	
	However, for twisted derived categories, this is not the correct construction to consider in general. 
		Let $A$ and $E$ be curved algebras.
	We state the following theorem.
	
	\begin{thm}\label{T:twbim}
	 There is a model structure on $\DGMod[A]{E}$ where a morphism $f:M\to N$ is 
		\begin{itemize}
			\item[(1)] a \emph{weak equivalence}  if
			$$f_\ast:\Homdg_{A^\op \otimes E}(T,M)\to \Homdg_{A^\op \otimes E}(T,N)$$
			is a quasi-isomorphism for all $T\in \Tw(A^\op)\otimes \Tw(E)$;
			\item[(2)] a \emph{fibration} if it is surjective;
			\item[(3)] a \emph{cofibration} if it has the left lifting property with respect to acyclic fibrations.
 		\end{itemize}
		We will denote the homotopy category by $\DII(A,E)$ and denote hom spaces $\Homdg_{A,E}$.
	\end{thm} 
\begin{proof}
	There is an equivalence of abelian dg categories
	\begin{align}
		\Phi:\DGMod[A]{E}&\to \DGMod{\Tw(A^\op)\otimes \Tw(E)}\label{Phi}\\
		M&\mapsto \Homdg_{A^\op \otimes E}(-,M) \label{Phi functor}.\notag
	\end{align}
	The model structure defined above is obtained by equipping the category of $\Tw(A^\op)\otimes \Tw(E)$-modules with the standard projective model structure and then equipping $\DGMod[A]{E}$ with the transported model structure along $\Phi$.
\end{proof}

	We denote this model structure by $\DGModII[A]{E}$ and denote the derived hom spaces by $\RHom_{A,E}$.
	
	We emphasize again that this is different in general from $\DGModII{A^\op\otimes E}$. 
	The reason is that not every twisted module over $A^\op \otimes E$ arises from tensor products of twisted $A$-modules and twisted $E$-modules, see Example \ref{eg:counterexample to condition 1}.
	
\begin{rem}\label{rem:bousfield}
By definition $\DGModII[A]{E}$ is a right Bousfield localization of $\DGModII{A^\op\otimes E}$.
\end{rem}	
\begin{prop}
	\begin{enumerate}
		\item Let $M\in \Tw(A^\op)$ and $N\in \Tw(E)$. Then $M\otimes N$ is cofibrant with respect to the model structure defined in Theorem \ref{T:twbim}. All cofibrant modules are retracts of filtered colimits of such tensor products taken along cofibrations. 
		\item\label{it:coacyclic} Any injection $f: M \to N$ in $\DGModII[A]{E}$ with cofibrant cokernel is a cofibration and vice versa.
	\end{enumerate}
\end{prop}
\begin{proof}
	\begin{enumerate}
		\item The cofibrant objects in $\DGModII[A]{E}$ are obtained by taking the image of $\Cof(\Tw(A^\op)\otimes \Tw(E))$ under $\Phi$ and that proves the first assertion.
		\item The cofibrations in $\DGModI{\Tw(A^\op)\otimes \Tw(E)}$ are exactly the injections with cofibrant cokernel thus the cofibrations in $\DGModII[A]{E}$ can be characterized similarly. \qedhere
	\end{enumerate}
\end{proof}
\begin{prop}\label{P:compactbimodules}
	The category $\DII(A, E)$ is compactly generated by $\Ho(\Tw(A^\op)\otimes \Tw(E))$.
\end{prop}
\begin{proof}
	The category $\D(\Tw(A^\op)\otimes \Tw(E))$ is compactly generated by the representable $\Tw(A^\op)\otimes \Tw(E)$-modules and by the proof of Theorem \ref{T:twbim} the representables correspond exactly to the objects in the subcategory $\Tw(A^\op)\otimes \Tw(E)\subset \DII(A, E)$ under $\Phi$.
\end{proof}
	\begin{cor}\label{comp:bimod}
		There is an equivalence $\DII(A,E)\simeq \D(\PerfII(A^\op)\otimes \PerfII(E))$ of triangulated categories
	\end{cor}
	\begin{proof}
		The functor $\Phi$, see Equation \ref{Phi}, induces an equivalence $\DII(A,E)\simeq \D(\Tw(A^\op)\otimes \Tw(E))$ of triangulated categories. The inclusion $\Tw(A)\subseteq \PerfII(A)$ is a derived Morita equivalence of dg categories thus $\D(\Tw(A^\op)\otimes \Tw(E))\simeq \D(\PerfII(A^\op)\otimes \PerfII(E))$.
	\end{proof}
Our next task is to generalize the Koszul duality statements, \cref{T:KdCobar} and \cref{T:KdBar}, to bimodules.
We treat the case of the cobar construction first.

	\subsection{Bimodule Koszul duality for the cobar construction}\label{sect:cobar}
	
		Throughout this section $C$ and $D$ will denote two pc curved algebras, 
	and the notation $\DGMod[C]{D}$ will be understood to mean pc bimodules, i.e.\ pc left $C$-modules that are also pc right $D$-modules.
	
	We now obtain an analogue of \cref{P:resolutiononesided} for pc bimodules by defining the following adjunction: 
	\[
	G \colon (\DGMod[C]{D})^\cop \rightleftarrows \DGMod[\Omega C]{\Omega D} \lon F.
	\]
	Let $\xi_C \in \MC(C^\op \otimes \Omega C^\op)$ 
	and $\xi_{D} \in \MC(D \otimes \Omega D)$ be the canonical Maurer--Cartan elements corresponding to the counits $\hB \Omega C^\op \to C^\op$ and $\hB \Omega D \to D$ of the adjunction $\Omega \dashv \hB $. 
	Define
	\[
	\xi := \xi_C \otimes 1 + 1 \otimes \xi_{D} \in C^\op \otimes \Omega C^\op \otimes D \otimes \Omega D,
	\]
	then $\xi \in \MC(C^\op \otimes \Omega C^\op \otimes D \otimes \Omega D)$.
	The functor $G$ sends a pc dg $C$-$D$-bimodule $N$ to
	\[
	GN \coloneqq ( \Omega C \otimes N^* \otimes\Omega D)^{[\xi]},
	\]
	where as before, $GN$ is a $(\Omega C,\Omega D)$-bimodule. 
	It is the bimodule cobar-construction of the $C^\op \otimes D$-module $N$. It can be written as the direct sum totalization of a double complex as follows (noting that the arrows are pre-differentials which need not to square to 0):
	\begin{equation}\label{doublecomplex}
		\xymatrix{
			\ar[d]&\ar[d]&\ar[d]\\
			N^* \otimes \widebar{D}^2\ar[d]
			& \widebar{C}\otimes N^*\otimes \widebar{D}^2\ar[d]\ar[l]
			&\widebar{C}^2\otimes N^*\otimes \widebar{D}^2\ar[d]\ar[l]
			&\ar[l]\cdots\\
			N^*\otimes \widebar{D}\ar[d]&\widebar{C}\otimes N^*\otimes \widebar{D}\ar[d]\ar[l]&\widebar{C}^2\otimes N^*\otimes \widebar{D}\ar[d]\ar[l]&\ar[l]\cdots\\
			N^*&\widebar{C} \otimes N^*\ar[l]&\widebar{C}^2\otimes N^*\ar[l]&\ar[l]\cdots\\
		}
	\end{equation}
	This is, in fact, a triple complex where the third pre-differential (not indicated explicitly) is induced on each term by the internal pre-differentials in $C, D$ and $N$.  
	The $n$th row of the above complex is $G_C(N\otimes \widebar{D}^n)$ and the $n$th column is $G_{D}(\widebar{C}^n\otimes N)$. 
	
	Similarly, the functor $F$ sends a curved $(\Omega C,\Omega D)$-bimodule $M$ to the pc dg $(C,D)$-bimodule
	\[
	FM \coloneqq (C \otimes M^* \otimes D)^{[\xi]},
	\]
	and the functors $F$ and $G$ are adjoint in the sense that there is a natural isomorphism of dg vector spaces
	\begin{align}
		\Hom_{\DGMod[\Omega C]{\Omega D}}(GN, M) \cong \Hom_{\DGMod[C]{D}}(FM, N)\label{adjunction iso}		
	\end{align}
	Indeed, both sides above are identified with the dg vector spaces $(N\otimes M)^{[\xi]}$.
	
	Composing the functors, we obtain the $(C,D)$-bimodule 
	\[
	FGN = (C \otimes \Omega C^\ast \otimes N \otimes \Omega D^\ast \otimes  D)^{[\xi\otimes 1+1\otimes \xi]}.
	\]
	This is a double complex obtained from (\ref{doublecomplex}) by tensoring each entry with $C$ from the left and $D$ from the right.
	
	Furthermore, this new double complex can be `augmented' by adding to it as a $(-1)$-st row the cobar-resolution $F_CG_C(N)$ of the $C$-module $N$. 
	We obtain 
	\begin{equation*}\label{eq:doublecomplex}
		\xymatrix{
			&\ar[d]&\ar[d]&\ar[d]\\
			&C\otimes N\otimes \widebar{D}^2 \otimes D\ar[d]
			&C \otimes \widebar{C}\otimes N \otimes \widebar{D}^2 \otimes D\ar[d]\ar[l]
			&C \otimes \widebar{C}^2\otimes N\otimes \widebar{D}^2 \otimes  D\ar[d]\ar[l]
			&\ar[l]\cdots\\
			&C \otimes N\otimes \widebar{D}\otimes D\ar[d]
			&C \otimes \widebar{C}\otimes N\otimes \widebar{D} \otimes D
			\ar[d]\ar[l]
			&C\otimes \widebar{C}^2\otimes N\otimes \widebar{D} \otimes D
			\ar[d]\ar[l]&\ar[l]\cdots\\
			&C \otimes N \otimes D \ar[d]
			&C \otimes \widebar{C} \otimes N\otimes D\ar[l]\ar[d]
			&C \otimes \widebar{C}^2\otimes N\otimes D\ar[l]\ar[d]
			&\ar[l]\cdots\\
			&		C\otimes N 
			&C \otimes \widebar{C} \otimes N \ar[l]
			&C \otimes \widebar{C}^2\otimes N\ar[l]
			&\ar[l]\cdots\\
		}
	\end{equation*}
	
	The resulting total complex can be viewed as the cofiber of a map $FGN \to F_CG_C(N)$.
	We canonically truncate in the vertical and horizontal direction and as in the one-sided case we obtain an inverse system of bounded acyclic complexes of finite length.  This shows that this cofiber is contra-acyclic as a $(C,D)$-bimodule. Note that the truncations are indeed complexes, see the discussion in the proof of Theorem \ref{T:KdCobar}.
	
	Thus, we conclude that the following result holds:
	
	\begin{prop}\label{P:resolutiontwosided}
		The adjunction unit $FGN \to N$ is a weak equivalence of curved $(C,D)$-bimodules.
	\end{prop}
	
	One proves similarly: 
	\begin{prop}\label{P:resolutiontwosidedunit}
		The adjunction counit $GFM \to M$ is a weak equivalence of $\Omega C^\op \otimes \Omega D$-modules.
	\end{prop}
	\begin{proof}
		The argument proving Proposition \ref{P:resolutiontwosided} proves that the unit map has a contra-acyclic cone and Proposition \ref{prop:useful}.\ref{it:coacyclic} shows this is a weak equivalence.
	\end{proof}
	Note that the proposition is equally true with the same proof if we consider the model structure on $(\Omega C, \Omega D)$-bimodules.
	
	We can now formulate the following bimodule version of \cref{T:KdCobar}. We prove here the version for $\Omega C^\op \otimes \Omega D$-modules, we will see in Corollary \ref{C:KdCobarAgain} that the same adjunction holds for $(\Omega C,\Omega D)$-bimodules.

	\begin{thm}\label{T:bimodKdCobar}
		The functor $G$ is left adjoint to $F$ and they form a Quillen equivalence
		\[
		G \colon (\DGMod[C]{D})^\cop \rightleftarrows \DGModII{ \Omega C^\op \otimes \Omega D} \lon F.
		\] 
	\end{thm}
	
	\begin{proof}
		We have an adjoint pair of functors $(G,F)$ between $C^{\op}\otimes D$-modules and $\Omega C^\op \otimes \Omega D$-modules; it has already been argued above that this is indeed an adjoint pair. 
		Moreover, the functor $G$ preserves cofibrations. First observe that any surjection of pc modules is sent to an injection.
		Next, any pc module is a limit of finite-dimensional modules and thus sent by $G$ to a colimit of twisted complexes over $\Omega C^\op \otimes \Omega D$.
		In particular, the image of the kernel, which is the cokernel of the image, is always cofibrant.
		
	Suppose $f$ is an acyclic cofibration, then $\operatorname{cone}(f)$ is coacyclic (note, contraacyclic objects become coacyclic in the opposite category). The functor ${G}$ takes absolutely acyclic $(C,D)$-bimodules to absolutely acyclic $\Omega C^\op \otimes \Omega D$-modules, which by Proposition \ref{prop:useful} are weakly trivial. 
	Furthermore, ${G}$ preserves arbitrary coproducts and a coproduct of weakly trivial $\Omega C^\op \otimes \Omega D$-modules is again weakly trivial. 
		
		
		This shows that $(G,F)$ is a Quillen adjunction. 
		By Propositions \ref{P:resolutiontwosided} and \ref{P:resolutiontwosidedunit} $F$ and $G$ induce an isomorphism at the level of homotopy categories. 
	\end{proof}

	\begin{cor}\label{T:monoidalCobar}
		For any pc curved algebras $C$ and $D$, there is an equivalence of triangulated categories between 
		$\DII(\Omega C^\op \otimes \Omega D)$ and $\DII(\Omega (C^\op \otimes D))$.
	\end{cor}
	
	\begin{proof}
		By Koszul duality there is a Quillen equivalence 
		\[
		(\DGMod[C]{D})^\cop \rightleftarrows \DGModII{\Omega (C^\op \otimes D)},
		\]
		so by \cref{T:bimodKdCobar} there is an equivalence of homotopy categories $\Ho (\DGMod[\Omega C]{\Omega D})$ and $\Ho (\DGMod{\Omega (C^\op \otimes D)})$.
	\end{proof}
	\begin{rem}
		Note that this subsection would have simplified drastically if the functor $\Omega$ was quasi-strong monoidal. This is known for local pc algebras and those dual to pointed coalgebras \cite{holstein2025enriched}, but not in general.
	\end{rem}
	\subsection{Bimodule Koszul duality for the bar constructions}\label{sect:bimkoszul}
	
	Throughout this section, $A$ and $E$ will denote two curved algebras, and $\DGMod[A]{E}$ will denote the category of curved $(A,E)$-bimodules.

 We define two functors
	\[
	\widecheck{G} \colon (\DGMod[\cB A]{\cB E})^\cop \rightleftarrows \DGMod[A]{E} \lon \widecheck{F}.
	\]
	Let $\xi_A \in \MC(A^\op \otimes \cB A^\op)$ and $\xi_{E} \in \MC(E \otimes \cB E)$ be the canonical Maurer-Cartan elements corresponding to the counits $\Omega \cB A^\op \to A^\op$ and $\Omega \cB E \to E$ of the adjunction $\Omega \dashv \cB $. 
	Define
	\[
	\xi := \xi_A \otimes 1 + 1 \otimes \xi_{E} \in A^\op \otimes \cB A^\op \otimes E \otimes \hB E,
	\]
	then $\xi \in \MC(A^\op \otimes E \otimes \cB A^\op \otimes \cB E)$.
	The functor $\widecheck{F}$ sends a curved $(A,E)$-bimodule $M$ to
	\[
	\widecheck{F}M \coloneqq (\cB A \otimes M^*  \otimes \cB E)^{[\xi]}
	\]
	and the functor $\widecheck{G}$ sends a pc curved $(\cB A,\cB E)$-bimodule $N$ to
	\[
	\widecheck{G}N \coloneqq (A \otimes N^* \otimes E)^{[\xi]}.
	\]

			Our goal now is to prove the following theorem, which is the bimodule version of \cref{T:KdBar}.
	\begin{thm}\label{T:quillenbimodule}
		There is a Quillen equivalence
		$$\widecheck{G}\colon (\DGMod[\cB A]{\cB E})^\cop \leftrightarrows \DGModII[A]{E}\colon \widecheck{F}.$$
	\end{thm}
\begin{rem}
If $A$ and $E$ are dg algebras,  we can similarly define two functors
\[
\widehat{G} \colon (\DGMod[\hB A]{\hB E})^\cop \rightleftarrows \DGMod[A]{E} \lon \widehat{F}
\]
exactly as above, except replacing every occurrence of $\cB $ by $\hB $. In this case, an analogue of Theorem \ref{T:quillenbimodule} holds  with $\DGModII[A]{E}$ replaced with the category of $(A,E)$-bimodules with the ordinary projective model structure (so that weak equivalences are ordinary quasi-isomorphisms). We do not need this result and so will not give a proof; however note that in the case $A=E$ and $A$ is augmented, it is proved in \cite[Proposition 2.4]{kel21a}.
\end{rem}
		The proof relies on a couple of preliminary results and will be take up most of the rest of this section.
	
	Let $\cat{pccuAlg}$ denote the category of pc curved algebras and $V$ be a finite-dimensional graded vector space.
	\begin{prop}\label{MC-elements of bimodules}
		There is a bijection
		\begin{align*}
			&\Hom_{\cat{pccuAlg}}(\cB A^\op\hat{\otimes} \cB E,\End(V))\\
			&\cong \{(\alpha_1,\alpha_2)\in \MC(A\otimes \End(V))\times \MC(\End(V)\otimes E)\ |\ [\alpha_1\otimes 1_E,1_A\otimes \alpha_2]=0\}.
		\end{align*}
	\end{prop} 
	\begin{proof}
	Let
	\begin{align}
		\textbf{M}:=&\{(f,g)\in \Hom_{\cat{pccuAlg}}(\cB{A}^\op,\End(V))\times \Hom_{\cat{pccuAlg}}(\cB{E},\End(V))\notag \\
		&|f(x)g(y)-g(y)f(x)=0, \text{for all }x\otimes y\in \cB{A}\hat{\otimes} \cB{E} \}, \label{the set M}
	\end{align}
	then there is a bijection
		\begin{align*}
			\textbf{M}&\cong \Hom_{\cat{pccuAlg}}(\cB{A}^\op\hat{\otimes} \cB{E},\End(V)),\\
			(f,g)&\mapsto [x\otimes y\mapsto f(x)g(y)].
		\end{align*}
		By \cite[Proposition 2.6]{gl21}, every $(f,g)\in \textbf{M}$ corresponds to a pair $(\alpha_1,\alpha_2)\in \MC(A\otimes \End(V))\times \MC(\End(V)\otimes E)$ and the condition that $f$ and $g$ commute, as described in \ref{the set M}, translates to the condition that $[\alpha_1\otimes 1_E,1_A\otimes \alpha_2]=0$. 
	\end{proof}
	\begin{lem}\label{equivalent MC}
		Let $\alpha_1\in \MC(A\otimes \End(V))$ and $\alpha_2\in \MC(\End(V)\otimes E)$. 
		The following are equivalent:
		\begin{itemize}
			\item[(1)] $\alpha_1\otimes 1_E+1_A\otimes \alpha_2\in \MC(A\otimes \End(V)\otimes E)$ and
			\item[(2)] $[\alpha_1\otimes 1_E,1_A\otimes \alpha_2]=0$
		\end{itemize}
	\end{lem}
	\begin{proof}
		This is just a matter of unravelling the Maurer-Cartan equation;
		\begin{align*}
			&(\alpha_1\otimes 1_E+1_A\otimes \alpha_2)^2+d_{A\otimes \End(V)\otimes E}(\alpha_1\otimes 1_E+1_A\otimes \alpha_2)\\
			&=\alpha_1^2\otimes 1_E+d_{A\otimes \End(V)}(\alpha_1)\otimes 1_E+1_A\otimes d_{\End(V)\otimes E}(\alpha_2) +1_A \otimes \alpha_2^2
			+[\alpha_1\otimes 1_E,1_A\otimes \alpha_2]\\
			&=-(h_A\otimes \operatorname{id}_{\End(V)\otimes E}+\operatorname{id}_{A\otimes \End(V)}\otimes h_E) +[\alpha_1\otimes 1_E,1_A\otimes \alpha_2].
		\end{align*}
		From the computation above, it is clear that the vanishing of $[\alpha_1\otimes 1_E,1_A\otimes \alpha_2]$ is equivalent to $\alpha_1\otimes 1_E+1_A\otimes \alpha_2\in \MC(A\otimes \End(V)\otimes E)$.
	\end{proof}
	\begin{defn}
		Let $A$ and $E$ be curved algebras and define
		\begin{align*}
			&\MC_\text{Bimod.}(A\otimes \End(V)\otimes E)\\
			&=\{(\alpha_1,\alpha_2)\in \MC(A\otimes \End(V))\times \MC(\End(V)\otimes E)\ |\ [\alpha_1\otimes 1_E,1_A\otimes \alpha_2]=0\}.
		\end{align*}
	\end{defn}
	Recall the following classical result from representation theory. See \cite[Lemma 3.4.2]{Etingof11} for a proof in the ungraded case that also works in the graded case.
	\begin{lem}\label{rep. lemma 1}
		Let $A$ be a graded algebra and let $V$ be a finite-dimensional graded $A$-module. There is a finite filtration
		\begin{align}
			0=V_0\subseteq V_1\subseteq ...\subseteq V_n=V,
		\end{align}
		such that $V_{i+1}/V_i$ is a simple $A$-module. 
		
		Furthermore, under the isomorphism of vector spaces $V\cong \bigoplus^n_{i=0} V_{i+1}/V_i$ the $A$-module $V$ can be expressed as an iterated extension of the simple $A$-modules $V_{i+1}/V_i$.
	\end{lem}
		We will also need the following lemmas from graded algebra which are the reason for our restrictions on the ground field. Recall that we are always assuming 	that $k$ is a perfect field and that $\operatorname{char}(k) \neq 2$ if the grading group is $\Z/2$.
		
		\begin{lem}\label{lem:sstensor}
			Let $A$ and $E$ be finite-dimensional graded semi-simple algebras. Then $A \otimes E$ is graded semi-simple.
					\end{lem}
					\begin{proof}
					By the graded Artin-Wedderburn theorem \cite[Theorem 2.10.10]{NastasescuVanOystaeyen2004} we may assume  $A$ and $E$ are (graded) matrix algebras over division rings.
												
						We now assume the grading group is $\Z$. 
						As the algebras are finite-dimensional, the division rings must be concentrated in degree 0.
					Thus, $A \otimes E$ is a matrix algebra over the tensor product of division rings, and the latter is itself semi-simple as a product of separable algebras (using that $k$ is perfect)  \cite[\S 10.5, \S 10.7]{Pierce1982}.
												As a graded matrix algebra over a semi-simple algebra, $A \otimes E$ is graded semi-simple.
												
%
%
				Let us now assume we have a $\Z/2$-grading. Let $K, L$ be the even centers of $A$ and $E$. 
						Then the tensor product of separable field extensions $K \otimes_k L$ is isomorphic to a product of fields $\prod_i F_i$.
						It follows that we may write $A \otimes_k E$ as $\prod_i (A\otimes_K F_i) \otimes_{F_i} (E \otimes_L F_i)$, which is a product of graded tensor products of central simple algebras. The tensor product of $\Z/2$-graded central simple algebras is graded simple by a theorem of Wall \cite[Theorem 2]{Wall1964}. Thus $A\otimes E$ is graded semi-simple.					
					\end{proof}
	\begin{lem}\label{lem:separablealgebras}
		Let $A$ and $E$ be finite-dimensional graded algebra. Then any simple graded $A\otimes E$-module $L$ is a direct summand of a tensor product of a simple graded $A$-module and a simple graded $E$-module.
	\end{lem}
	\begin{proof}
		Let $ J_A=\operatorname{rad}(A)$ and $J_E=\operatorname{rad}(E)$ be the graded Jacobson radicals. As in the ungraded case $J_A, J_E$ are nilpotent and $A/J_A$ and $E/J_E$ are (graded) semi-simple.
		We consider the ideal $I=J_A\otimes E+A\otimes J_E \subseteq A\otimes B $ which is nilpotent as the radicals are and homogeneous.  
		Thus $I$ acts trivially on the simple module $L$ and $L$ is a module over the graded algebra $R/I$ where $R = A\otimes B$.
		
	By construction $R/I \cong (A/J_A)\otimes_k(E/J_E)$ and by Lemma \ref{lem:sstensor} the right hand side is graded semi-simple.
		Thus $L$ is a simple graded module over a graded semi-simple algebra and up to shift it is a direct summand of $A/J_A \otimes E/J_E$, proving the lemma.
%
%
	\end{proof}
	Whenever $\cat{A}$ is a dg category, we write $\widehat{\tw}(\cat{A})$ for the (strict) idempotent completion of $\tw(\cat{A})$.
	\begin{prop}\label{simple twisted modules}
		Let $T=(A\otimes V\otimes E,d_{A\otimes V\otimes E}+\alpha)\in \Tw(A^\op\otimes E)$ with $\alpha\in \MC(A\otimes \End(V)\otimes E)$. 
		\begin{enumerate}
			\item[(1)] If $\alpha\in \MC_\text{Bimod.}(A\otimes \End(V)\otimes E)$ then $T\in \widehat{\tw}(\Tw(A^\op)\otimes \Tw(E))$.			\item[(2)]  If
			$T\in \Tw(A^\op)\otimes \Tw(E)$ then $\alpha\in \MC_\text{Bimod.}(A\otimes \End(V)\otimes E)$. 
		\end{enumerate}
	\end{prop}
	\begin{proof}
		\begin{enumerate}
			\item[(1)] By assumption, $\alpha=(\alpha_1,\alpha_2)$ with $\alpha_1=\sum^m_i a_i\otimes E_i\in A\otimes \End(V)$ and $\alpha_2=\sum^n_j F_j\otimes e_j\in \End(V)\otimes E$. Let $X, Y\subseteq \End(V)$ denote the (graded) algebras generated by $\{E_i\},\{F_j\}$, respectively. As the generators of $X$ and $Y$ commute, there is a well-defined map
			\begin{align*}
				X\otimes Y\to \End(V)\\
				x\otimes y\mapsto xy
			\end{align*}
			making $V$ a $X\otimes Y$-module.
			
			Now let $0=V_0\subseteq V_1\subseteq ...\subseteq V_k=V$ be a filtration of $X\otimes Y$-invariant subspaces of $V$.
			It follows from Lemma \ref{rep. lemma 1}, 
			that $V$ can be expressed as an iterated extension of simple $X\otimes Y$-modules $W_l=V_l/V_{l-1}$. As we are working over a perfect field, Lemma \ref{lem:separablealgebras} applies, hence $W_l$ must be a summand of some $X\otimes Y$-module $M_l\otimes N_l$, where $M_l$ is a simple $X$-module and $N_l$ is a simple $Y$-module.
			\\Assume that $W_l\cong M_l\otimes N_l$ for all $0\leq l\leq k$, then $T$ is an iterated extension of finitely generated twisted $(A,E)$-bimodules of the form 
			\begin{align}
				\left(A\otimes M_l\otimes N_l\otimes E,d_{A\otimes M_l\otimes N_l\otimes E}+\sum^m_{i=1} a_i\otimes E^l_i\otimes \operatorname{id}_{N_l}\otimes 1_E+\sum^n_{j=1} 1_A\otimes \operatorname{id}_{M_l}\otimes F^l_j\otimes e_j\right).\label{twist:simple form}
			\end{align}
			If $W_l$ is a proper summand for some $0\leq l\leq k$, the $(A,E)$-bimodule $(A\otimes W_l\otimes E,d)$ obtained by the procedure above, must be a summand of some $(A,E)$-bimodule of the form \ref{twist:simple form}, which is clearly an object in $\Tw(A^\op)\otimes \Tw(E)$. Therefore, $(A\otimes W_l\otimes E,d)\in \widehat{\tw}(\Tw(A^\op)\otimes \Tw(E))$ which in turn implies that $T\in \widehat{\tw}(\Tw(A^\op)\otimes \Tw(E))$. That proves the first assertion.
			\item[(2)] If $T\in \Tw(A^\op)\otimes \Tw(E)$, then by Lemma \ref{equivalent MC}, $\alpha\in \MC_\text{Bimod.}(A\otimes \End(V)\otimes E)$.
		\end{enumerate}
	\end{proof}
	Recall that a pc curved $C$-module structure on a vector space $V$ corresponds to a map of curved pc algebras $C\to \End(V)$. 
	\begin{cor}\label{MC-elements and twisted modules}
		\begin{enumerate}
			\item Let $L$ be a finite-dimensional pc curved $(\cB{A},\cB{E})$-bimodule, then $\widecheck{G}(L)\in \widehat{\tw}(\Tw(A^\op)\otimes \Tw(E))$.
			\item Let $T\in \Tw(A^\op)\otimes \Tw(E)$, then there exists a finite-dimensional  $(\cB{A},\cB{E})$-bimodule $K$ such that $T=\widecheck{G}(K)$.
		\end{enumerate}
	\end{cor}
	\begin{proof}
		\begin{enumerate}
			\item Let $L$ be a finite-dimensional pc curved $(\cB{A},\cB{E})$-bimodule, then by Proposition \ref{MC-elements of bimodules} and Lemma \ref{equivalent MC} we obtain a Maurer-Cartan element $\alpha_L\in \MC_\text{Bimod.}$, which by Proposition \ref{simple twisted modules} corresponds to a $(A,E)$-bimodule $T\in \widehat{\tw}(\Tw(A^\op)\otimes \Tw(E))$, where $T=\widecheck{G}(K)$.
			\item Given a $T\in \Tw(A^\op)\otimes \Tw(E)$ it follows from Proposition \ref{simple twisted modules} that there is a corresponding $\alpha \in \MC_\text{Bimod.}(A\otimes \End(V)\otimes E)$ for some finite-dimensional graded vector space $V$. Under the one-to-one correspondence established in Proposition \ref{MC-elements of bimodules} we obtain a map of pc curved algebra
			$\cB A^\op\hat{\otimes}\cB E\to \End(V)$. It is a standard result, that such a map gives rise to a finite-dimensional pc curved $(\cB A,\cB E)$-bimodule $K$, 
			see Proposition \ref{simple twisted modules} (2). Tracing through these correspondences, we conclude that $	\widecheck{G}(K)=(A\otimes K^\ast\otimes E)^{[\xi]}\in \Tw(A^\op)\otimes \Tw(E)$.
		\end{enumerate}
	\end{proof}
Now we introduce the bimodule of 1-forms in the curved setting, which will be used in the proof that the unit of the adjunction $\check G \dashv \check F$ is a weak equivalence, see Proposition \ref{adj unit:weak eq}.

Let us recall that for a dg algebra $A$, its  bimodule of noncommutative 1-forms $\Omega_A$ is defined as the kernel of the multiplication map $A\otimes A\to A$ so there is a short exact sequence of dg $A$-bimodules
\begin{equation}\label{eq:noncommdiff}\xymatrix{
		0\ar[r]&\Omega_A\ar^i[r]&A\otimes A\ar[r]&A\ar[r]&0.
	}
\end{equation}
The image of $i$ inside $A\otimes A$ is generated by elements of the form $a\otimes 1-1\otimes a$ with $a\in A$. Set $\eta(a):=i^{-1}(a)\in \Omega_A$; the map $A\to \Omega A:a\mapsto \eta(a)$ is known as the universal derivation of $A$.

Denote the total complex of $\Omega_A\to A\otimes A$ by $W_A$; thus $W_A\cong\Omega_A[1]\oplus A\otimes A$ with the total differential induced by those on $\Omega_A,A\otimes A$ and $i$. It is clear that $W_A$ is quasi-isomorphic, as an $A$-bimodule, to $A$.

Let us now assume that $A$ is a (ordinary or pc) curved algebra with the curvature element $c$, i.e. for $a\in A$ we have $d^2(a)=[c,a]$ and $d(c)=0$. In that case, we can still form $\Omega_A$, and both $\Omega_A$ and $A\otimes A$ will have compatible differentials (possibly not squaring to zero). Furthermore, $A$ is still a (curved) $A$-bimodule, but neither $\Omega_A$ nor $A\otimes A$ are curved $A$-bimodules. The sequence (\ref{eq:noncommdiff}) is \emph{not} a sequence of curved $A$-bimodules. Nevertheless, we will show that the differential in $W_A$ can be deformed so that it becomes a curved $A$-bimodule.

Namely, define the total differential $D_c$ on the graded space underlying $W_A$ as follows. For $\eta(a)\in\Omega A[-1]$ set $D_c(\eta(a))=-\eta(da)+a\otimes1-1\otimes a$ (so this piece of the differential is the same as in the dg case). For $1\otimes 1\in A\otimes A$ set $D_c(1\otimes 1)=\eta(c)\in\Omega_A[1]$. Since the elements $\eta(a)$ and $1\otimes 1$ generate $W_A$ as a graded $A$-bimodule, this defines $D_c$ on the whole of $W_A$. Schematically, the curved structure on $W_A$ can be represented as the following diagram:

\[
\begin{tikzpicture}[>=stealth, baseline=(current bounding box.center)]
	\node (L) at (0,0) {$\Omega_A$};
	\node (R) at (4,0) {$A\otimes A$};
	
	\draw[->] (L) -- node[above] {$i$} (R);
	
	\draw[->] (R) to[out=150,in=30] node[above] {$h_c$}  (L);
	
	\draw[->]
	(L) .. controls +(-1,0.8) and +(-1,-0.8) ..
	node[left] {$d$} (L);
	
	\draw[->]
	(R) .. controls +(1,0.8) and +(1,-0.8) ..
	node[right] {$d$} (R);
\end{tikzpicture}
\]
In other words, the total differential 	in $W_A$ is modified from the dg case by adding the arrow  $h_c$ from $A\otimes A$ to $\Omega_A$ associated to the curvature $c$. Then the following result holds.
\begin{lem}\label{lem:curved}
	Given a curved algebra (pc or discrete) $A=(A,d,c)$, the graded vector space $W_A$ supplied with the differential $D_c$ as above, is a curved $A$-bimodule.	Moreover, the surjective map $W_A\to A$ indiced by the multiplication map $A\otimes A\to A$, has contractible kernel. 
\end{lem}
\begin{proof}
	To show that $(W,D_c)$ is a curved $A$-bimodule we should prove that for any $x\in W_A$, it holds that $D^2_c(x)=cx-xc$. Since $W_A$ is generated, as a graded $A$-bimodule, by elements of the form $\eta(a), a\in A$ and $1\otimes 1\in A\otimes A\subset W_A$, it suffices to check the required identity for $x=\eta(a)$ or $x=1\otimes 1$.
	
	For $x=\eta(a)$, we have
	\begin{align*}
		D^2_c(x)&=D_c(i(x)-d(x))\\&=D_c(a\otimes 1-1\otimes a-\eta(da))\\
		&=[da\otimes 1-1\otimes da+a\eta(c)-\eta(c)a]-[i\eta(da))-d\eta(da)]\\
		&=[da\otimes 1-1\otimes da+a\eta(c)-\eta(c)a]-[da\otimes 1-1\otimes da-\eta(d^2a)]=\\
		&=a\eta(c)-\eta(c)a+\eta(d^2a)\\
		&=[a,\eta(c)]+\eta[c,a]\\
		&=[c,\eta(a)]
	\end{align*}
	For $x=1\otimes 1$, we have
	\begin{align*}
		D^2_c(x)&=D_c(\eta(c))\\
		&=i(\eta(c))-d\eta(c)\\
		&=c\otimes 1-1\otimes c-\eta(d(c))=\\
		&=c\otimes 1-1\otimes c.
	\end{align*}
	as required.
	
	To show that $W_A$ is weakly equivalent to $A$ as a curved $A$-bimodule, note that there is an obvious map surjective map of curved $A$-bimodules $W_A\to A$ induced by the multiplication map $A\otimes A\to A$. The kernel of this map is isomorphic to the curved $A$-bimodule having the form
	
	\[
	\begin{tikzpicture}[>=stealth, baseline=(current bounding box.center)]
		\node (L) at (0,0) {$\Omega_A$};
		\node (R) at (4,0) {$\Omega_A$};
		
		\draw[->] (L) -- node[above] {$i$} (R);
		
		\draw[->] (R) to[out=150,in=30] node[above] {$h_c$}  (L);
		
		\draw[->]
		(L) .. controls +(-1,0.8) and +(-1,-0.8) ..
		node[left] {$d$} (L);
		
		\draw[->]
		(R) .. controls +(1,0.8) and +(1,-0.8) ..
		node[right] {$d$} (R);
	\end{tikzpicture}
	\]
	with $h_c$ and $i$ being mutually inverse isomorphisms. Such a curved $A$-bimodule clearly is contractible and so $W_A$ and $A$ are weakly equivalent.
\end{proof}
	\begin{prop}\label{adj unit:weak eq}
		The adjunction unit $\epsilon_N:\widecheck{F}\widecheck{G}(N)\rightarrow N$ is a weak equivalence of $(\cB A, \cB E)$-bimodules.
	\end{prop}
	\begin{proof}
		It follows from Lemma \ref{lem:curved} together with
		 \cite[Proposition 2.4]{gl21} that
		\begin{align*}
			\phi_A:(\cB A\otimes A^\ast\otimes \cB A)^{[\xi]}\to \cB A
		\end{align*}
		is a cofibrant resolution of $\cB {A}$ as $\cB {A}$-bimodule. Given a $(\cB {A},\cB {E})$-bimodule $N$, we get a weak equivalence of $(\cB {A},\cB {E})$-bimodules
		\begin{align*}
			\phi_A\otimes_A \id_N: (\cB {A}\otimes A^\ast\otimes N)^{[\xi]}\to N
		\end{align*}
		by tensoring $\id_N$ with $\phi_A$ over $\cB {A}$ on the left (tensoring preserves weak equivalences in the contraderived category). 
		Similarly, tensoring with $\phi_E$ on the right results in a weak equivalence
		\begin{align*}
			\phi_A\otimes_A \id_N\otimes_E \phi_E: (\cB A\otimes A^\ast\otimes N\otimes E^\ast\otimes \cB E)^{[\xi]}\to N.
		\end{align*}
		It can be verified that $\epsilon_N=\phi_A\otimes_A \id_N\otimes_E \phi_E$.
	\end{proof}
	\begin{prop}\label{adj counit:weak eq}
		The adjunction counit $\eta_M:\widecheck{G}\widecheck{F}(M)\to M$ is a weak equivalence of curved $(A,E)$-bimodules.
	\end{prop}
	\begin{proof}
		The morphism $\eta_M$ is a surjection, so the objective is to show that $\eta_M:\widecheck{G}\widecheck{F}(M)\to M$ is also a weak equivalence. Let $T\in \Tw(A^\op)\otimes \Tw(E)$, then $T=\widecheck{G}(V)$ for some finite-dimensional $V\in \DGMod[\cB A]{\cB E}$, see Corollary \ref{MC-elements and twisted modules}, and consider the commutative square:
		$$\begin{tikzcd}
			\Homdg_{A,E}(\widecheck{G}(V),\widecheck{G}\widecheck{F}(M))\ar[r,"(\eta_M)_\ast"]\ar[d, "\simeq"]& \Homdg_{A,E}(\widecheck{G}(V),M)\ar[d, "\simeq"]\\
			\Homdg_{\cB {A}, \cB {C}}(\widecheck{F}\widecheck{G}\widecheck{F}(M),V)& \Homdg_{\cB {A},\cB {C}}(\widecheck{F}(M),V)\ar[l, "\widecheck{F}(\eta_M)^\ast"']\ .
		\end{tikzcd}$$
		By the adjunction equations $\widecheck{F}(\eta_M)$ is a section of $\epsilon_{\widecheck{F}(M)}$ thus it must be a weak equivalence in $\DGMod[\cB A]{\cB E}$. As $\widecheck{F}(\eta_{M})$ is a weak equivalence of cofibrant objects and every object is fibrant, in particular $V$, $\widecheck{F}(\eta_{M})^\ast$ is a quasi-isomorphism.
		The 2-out-of-3 property of quasi-isomorphisms implies, that $(\eta_M)_\ast$ must be a quasi-isomorphism ensuring that $\eta$ is indeed a weak equivalence . 
	\end{proof}
		\begin{proof}[Proof of Theorem \ref{T:quillenbimodule}]
		The proof that $(\widecheck{G},\widecheck{F})$ is an adjunction, amounts to a computation identical to the one for $(G,F)$, see Section \ref{adjunction iso}.
		
		We now turn to the proof that the adjunction is Quillen.
		Suppose $f:M\to N$ is cofibration in $(\DGMod[\cB A]{\cB E})^\cop$, then we have to verify that $\widecheck{G}(f)$ is a cofibration, i.e. injective with cofibrant cokernel. Injectivity of $\widecheck{G}(f)$ is clear. For the second condition, first note that since $\widecheck{G}$ is a left adjoint, $\widecheck{G}(\operatorname{coker}(f))\cong \operatorname{coker}(\widecheck{G}(f))$. Every pc module can be expressed as a filtered limit of finite-dimensional pc modules, thus $\operatorname{coker}(f)=\varprojlim_{i\in I} P_i$, where $P_i$ is a finite-dimensional $(\cB A, \cB E)$-bimodule for all $i\in I$. Filtered limits become filtered colimits in the opposite category and since $\widecheck{G}$ is a left adjoint, we obtain $\widecheck{G}(\operatorname{colim}_{i\in I} P_i)=\operatorname{colim}_{i\in I} \widecheck{G}(P_i)$. It follows from Corollary \ref{MC-elements and twisted modules} (1), that $\widecheck{G}(P_i)\in \widehat{\tw}(\Tw(A^\op)\otimes \Tw(E))$, which
		implies that $\widecheck{G}(\operatorname{coker}(f))$ is cofibrant, see Theorem \ref{T:twbim}. 
		
		Suppose $f$ is an acyclic cofibration, then $\operatorname{cone}(f)$ is coacyclic (note, contraacyclic objects become coacyclic in the opposite category). The functor $\widecheck{G}$ takes absolutely acyclic $(\cB A,\cB E)$-bimodules to absolutely acyclic $(A,E)$-modules, which by Proposition \ref{prop:useful} are weakly trivial. Furthermore, $\widecheck{G}$ preserves arbitrary coproducts and a coproduct of weakly trivial $(A,E)$-bimodules is again weakly trivial. From these observations and the homotopy equivalence $\operatorname{coker}(\widecheck{G}(f))\simeq (A\otimes \operatorname{cone}(f)^\ast\otimes E)^{[\xi]}$,
		that $\widecheck{G}(f)$ is a acyclic cofibration. 
		 
	It was established, in Proposition \ref{adj unit:weak eq} and Proposition \ref{adj counit:weak eq}, that the adjunction morphisms 
		\begin{align*}
			\eta_M&:\widecheck{G}\widecheck{F}(M)\to M\\
			\epsilon_N&:\widecheck{F}\widecheck{G}(N)\rightarrow N
		\end{align*}
		are weak equivalences, which implies that the Quillen adjunction $(\widecheck{G},\widecheck{F})$ is a Quillen equivalence.
	\end{proof}
		\begin{rem}\label{rem:naivebarquillen}
		If we work with the twisted derived category of $A^\op \otimes E$-modules instead of the twisted derived category of bimodules there is still a Quilen adjunction
		\[
		\widecheck{G} \colon (\DGMod[\cB A]{\cB E})^\cop \rightleftarrows \DGModII{A^\op \otimes E} \lon \widecheck{F}.
		\]
		which is not in general an equivalence. In fact, this Quillen adjunction is a Quillen coreflection which factors through the Bousfield localization from $\DGModII{A^\op \otimes E}$ to $\DGModII[A]{E}$.
	\end{rem}
	
	\begin{cor}\label{C:bimoduleovercofibrant}
		There is a Quillen equivalence $\DGModII[A]{E} \simeq \DGModII{A' \otimes E'^\op}$ where $A' \to A$ and $E' \to E$ are cofibrant replacements in $\cuAlg$.
	\end{cor} 
	\begin{proof}
		There is a natural Quillen adjunction
		$\DGModII[A]{E} \to \DGModII[A']{E'}$
		induced by the cofibrant replacement map.
		By Theorem \ref{T:compactII} and Proposition \ref{prop:mcd2} this is a Quillen equivalence.
		
		It remains to compare $\DGModII[A']{E'}$ and $\DGModII{A'^\op \otimes E'}$.

		Let us first assume  $A' = \Omega C$ and $E' = \Omega D$.
		By Theorem \ref{T:quillenbimodule}
		$\DII(A',E')$ is then equivalent to $\Dctr{(\cB \Omega C^\op \otimes \cB \Omega D)}^\cop$.
		But by module Koszul duality and Corollary \ref{T:monoidalCobar} we have
		\[
		\Dctr{(\cB \Omega C^\op \otimes \cB \Omega D)}^\cop
		\cong \DII{\Omega(\cB \Omega C^\op \otimes \cB \Omega D)} \cong \DII(\Omega \cB \Omega C^\op \otimes \Omega \cB \Omega D)
		\]
		but as $\otimes$ preserves MC equivalences of pc algebras \cite[Theorem 12.11]{bl23},
		this is equivalent to $\DII(\Omega  C^\op \otimes \Omega D)$ which in turn is $\DII( (A')^\op \otimes E')$ as desired.
		
			If $A'$ and $E'$ are general cofibrant (of the second kind) algebras, we know they are homotopy equivalent to algebras of the form $\Omega C$ and $\Omega D$. 
		But then $A' \otimes E'$ is in fact homotopy equivalent to $\Omega C \otimes \Omega D$.
		It suffices to prove that two homotopic maps $f,g$ between any two algebras $A_1 \to A_2$ remain homotopic after tensoring with any algebra $B$. 
		A homotopy between $f$ and $g$ is a map $A_1 \to A_2 \otimes I$ where $I$ is any MC-interval (e.g. the 3-interval $I^3$ from \cite[Section 3.2]{bl23}) and tensoring it with $B$ we get the required homotopy between $f\otimes B$ and $g\otimes B$. 
		But homotopic algebras have equivalent twisted derived categories and this completes the proof.
		\end{proof}
	\begin{cor}\label{C:KdCobarAgain}
		There is a Quillen equivalence
		\[
		G \colon (\DGMod[C]{D})^\cop \rightleftarrows \DGModII[\Omega C]{\Omega D} \lon F.
		\] 
	\end{cor}

\begin{proof}
	This follows immediately from Corollary \ref{C:bimoduleovercofibrant} and Theorem \ref{T:KdCobar}.
\end{proof}
	
	\section{Hochschild cohomology of the  second kind}\label{sect:secondkind}
	In this section we define Hochschild cohomology of the  second kind, give an explicit formula and prove Morita invariance as well as compatibility with Koszul duality.

	We first recall classical Hochschild cohomology.
Let $\cat{A}$ be a dg category. 
The \emph{Hochschild cohomology complex} of $\cat{A}$ is $\HHc(\cat{A}) = \RHom_{\cat{A}\otimes\cat{A}^{\op}} (\cat{A},\cat{A})$.
This is a dg algebra (in fact a $B_\infty$-algebra) considered as the algebra of derived endomorphisms of $\cat{A}$ and an $\cat{A}$-bimodule. 
The \emph{Hochschild cohomology} of $\cat{A}$ is $\HH(\cat{A}) = H(\HHc(\cat{A}))$.

Hochschild cohomology is invariant under Morita equivalence, i.e.\
 the Yoneda embedding $\cat{A} \to \Perf(\cat{A})$ induces a quasi-isomorphism 
\[
\HHc(\cat{A}) \simeq \HHc(\Perf(\cat{A})).
\]

	\subsection{Hochschild cohomology of the second kind}
	Let now $A$ be a curved algebra.
		We will consider the category of $A$-bimodules with the twisted bimodule model structure from Theorem \ref{T:twbim}.
	
	From now on $\DGVect$ will be equipped with the projective model structure.
	\begin{lem}\label{L:quillen}
		The functor 
		$$\Homdg_{A, A}(-,-): (\DGModII[A]{A})^\op \otimes \DGModII[A]{A} \to \DGVect$$ 
		is a Quillen bifunctor.
	\end{lem}
	\begin{proof}
		By adjointness it suffices to show that 
		$$\otimes_k: \DGModII[A]{A} \otimes \DGVect \to \DGModII[A]{A}$$ is a Quillen bifunctor.
		So let $g: P \to Q$ be a cofibration in $\DGVect$ and $f: L \to M$ a cofibration in $\DGModII[A]{A}$ and consider 
		$$f \boxempty g: M \otimes P \amalg_{L \otimes P} L \otimes Q \to M \otimes Q.$$
		To check that $f \boxempty g$ is a cofibration if $f$ and $g$ are, it suffices to check the case that $g$ is a generating cofibration.
		
		Then we can write $g: k[n] \to \operatorname{cone}(\textsf{id}_{k[n]})$
		and we obtain $f \boxempty g: M \oplus L[1] \to M \oplus M[1]$ which is a cofibration if $f$ is as cofibrations in the model structure of Theorem \ref{T:twbim} are just injections.
		
		If $g$ is a generating acyclic cofibration it is of the form 
		$0 \to \operatorname{cone}(\textsf{id}_{k[n]})$. 
		Then $f \boxempty g: L \otimes Q \to M \otimes Q$ has cokernel $\operatorname{cone}(f)\otimes\operatorname{cone}(\textsf{id}_{k[n]})$ which is a totalization of a short exact sequence and thus coacyclic, and thus $f \boxempty g$ is a weak equivalence by Proposition \ref{prop:useful}.

		Finally, it is clear from the definition in Theorem \ref{T:twbim} that if $f$ is a weak equivalence so is $f \otimes \textsf{id}_B$ for $B$ bounded finitely generated  in $\DGVect$. Thus if $L \to M$ is a weak equivalence then $L \otimes P \to M \otimes P$ and its pushout along $L \otimes P \to L \otimes Q$ are weak equivalences, as is $L \otimes Q \to M \otimes Q$. By the 2-out-of-3 property, this gives a weak equivalence $f \boxempty g$.
	\end{proof}
	
	We now write $\RHomII_{A,A}(M,N)$ for the derived hom complex between $M$ and $N$ in $\DGModII[A]{A}$.
	\begin{defn}\label{D:main}
		Let $A$ be a curved algebra and $M$ a curved bimodule.
		The \emph{Hochschild cohomology complex of the  second kind} of $A$ with coefficients in $M$ is 
		\[
		\HHIIc(A, M) \simeq \RHomII_{A,A} (A,M). 
		\]
		We write $\HHIIc(A)$ for $\HHIIc(A,A)$ and $\HHII(A)$ for $H(\HHIIc(A,A))$.
	\end{defn}

		We are concerned with Hochschild cohomology in this paper, but also note the following definition:
		\begin{defn}
				Let $A$ be a curved algebra and $M$ a curved bimodule.
			The \emph{Hochschild homology complex of the  second kind} of $A$ with coefficients in $M$ is 
		$A \otimes^{L, \text{II}}_{A,A}M.$ where the tensor product is derived in the model category $\DGModII[A]{A}$.
		\end{defn}
		
		Here we choose the subscript $A,A$  of the tensor product to distinguish from the computation in $A\otimes A^\op$-modules that is not correct in general.
		
		The definition makes sense as the tensor product of $A$-bimodules is a left Quillen bifunctor into $\DGVect$ with the projective model by arguments similar to Lemma \ref{L:quillen}. 
		
		Next we construct an explicit complex computing Hochschild cohomology of the second kind.
	
	\begin{cor}\label{C:compute}
		$\HHII(A)$ is computed  as cohomology of the complex $(\cB A \otimes A)^{\xi}$ where the superscript induces two-sided twisting by $\xi_A$. 
	\end{cor}
	
	\begin{proof}
		The computation is the same as how one might compute Hochschild cohomology in terms of the usual tensor algebra.
		
		We have $\HHIIc(A) \cong \RHom_{\cB A \otimes \cB A^\op}(\cB A, \cB A)$ by Theorem \ref{T:quillenbimodule}
		and have to find a cofibrant replacement for $\cB A$ as an $A$-bimodule. Let us first consider the case when $A$ is augmented  (so the pc algebra $\cB A$ uncurved)
		We consider the short exact sequence
		\begin{equation}\label{eq:uncurved}
			\cB A \otimes \overline{A}^*[1] \otimes \cB A \to \cB A \otimes \cB A \xrightarrow{\mu} \cB A
		\end{equation}
		where $\mu$ is the multiplication map. Since the cone of the map $\mu$ is by construction the totalization of a short exact sequence
		this shows that $\cB A$
		is weakly equivalent in $\Dctr(\cB A \otimes \cB A)$ to $\cB A \otimes \overline{A}^*[1] \otimes \cB A \to \cB A \otimes \cB A$, which is moreover cofibrant as it is a topologically free pseudo-compact module.
		
		Thus we obtain 
		$$\HHIIc(A) \simeq (\Hom_k(A^*, \cB A), d) \simeq (A \otimes \cB A, d)$$ 
		and the induced differential $d$ is exactly the two-sided twisting by $\xi_A$, which is the usual Hochschild differential.
		
		Note that the totalization of the bicomplex $\cB A \otimes \overline{A}^*[1] \otimes \cB A \to \cB A \otimes \cB A$ is precisely $W_{\cB A}$. Thus, in the case when $A$ is not augmented (and so $\cB A$ is curved), $W_{\cB A}$ is cofibrant and, by Lemma \ref{lem:curved}  $W_{\cB A}$ is weakly equivalent to $\cB A$ as a pc $\cB A$-bimodule. The proof is finished just as in the uncurved case.
	\end{proof}
	
	\begin{rem}\label{remark:comparison}
		There is another definition of Hochschild cohomology of the second kind, cf. \cite{pp12} which is \emph{not} equivalent to our notion. To define it, let $\mathsf{B}'(A):=\bigoplus_{n=0}^{\infty} (\bar{A}^{\otimes n})^*[-1]$ 
		be the `semi-complete' bar-construction of $A$; it is a dg algebra that is neither discrete in general (because 
		$(\bar{A}^{\otimes n})^*$ is a pc vector space) nor pc (because an infinite direct sum of pc vector spaces is not pc).
		Nevertheless, the complex $(\mathsf{B}'(A)\hat{\otimes}A)^{\xi}$ makes sense and its cohomology can be taken as a definition of \emph{Hochschild cohomology of the second kind} (in the sense of Polishchuk-Positselski), which we wil denote by $\HH^{\mathrm{II}}_{\mathrm{PP}}(A)$.
		
		This complex only differs from the classical Hochschild cohomology complex in that the totalization of the natural bicomplex is a direct sum totalization rather than the direct product totalization.
		This is also sometimes called 
		\emph{compactly supported Hochschild cohomology}, 
		$\HH^*_{\mathrm{c}}(A)$, and the dual construction is called \emph{Borel-Moore Hochschild homology} $\HH_*^{\textrm{BM}}(A)$. 
		
		The pc bar-construction $\cB(A)$ is the pc completion of $\mathsf B'(A)$ and the ordinary bar-construction $\hB(A)$ is a further completion at the maximal ideal. It follows that there are maps of graded algebras 
		$$\HH^{\mathrm{II}}_{\mathrm{PP}}(A)\to \HHII(A)\to \HH(A).$$ In other words, $\HHII(A)$ is a kind of a half-way house between $\HH^{\mathrm{II}}_{\mathrm{PP}}(A)$ and  $\HH(A)$.
		We will come back to the comparison between the two constructions in Section \ref{sect:hhpp}. 
	\end{rem}	
	\begin{exmp}
		Let $A = \Lambda(e)$ be the dg algebra of dual numbers with a generator in degree 1 and trivial differential.	
		Let us also assume $k$ is algebraically closed to simplify the formulas.
		We may compute $\HHII(A)$ from the extended bar construction of $A$, which is the free pc algebra on one generator, and obtain (choosing notation in line with the Hochschild-Kostant-Rosenberg theorem)
		$$\HHII(A) \cong \prod_{\alpha \in k} k_\alpha[[\partial_e]]\otimes \Lambda(e)$$
		where $\partial_e$ is in degree 0.
		Comparing with $\HH(A)$ and $\HH^{\mathrm{II}}_{\mathrm{PP}}(A)$ as in
		Remark \ref{remark:comparison} we obtain the maps
		$$k[\partial_e] \otimes \Lambda(e) \to \prod_{k} k[[\partial_e]] \otimes \Lambda(e) \to k[[\partial_e]]\otimes \Lambda(e)$$
		Thus, in this example $\HHII$ is the pc completion of the Polishchuck-Positselski construction and the usual Hochschild cohomology is further completion around 0.
	\end{exmp}

\subsection{Morita invariance}

In this section we consider the natural dg enhancement of $\DII(A)$ provided by cofibrant objects in $\DGModII{A}$ and denote it by $\DIIdg(A)$.
	
	\begin{thm}\label{thm:hhperf2}
	For any curved algebra $A$ there are quasi-isomorphisms of dg algebras	$\HHIIc(A) \simeq \HHc(\PerfII(A)) \simeq \HHc(\DIIdg(A))$.
	\end{thm}
	Note that to make sense of the second quasi-isomorphism we should consider $\DIIdg(A)$ to be a small dg category, e.g.\ by assuming the existence of universes $\mathbb U \in \mathbb{W}$ and considering $A$-modules to be $\mathbb{U}$-small chain complexes with an action of $A$.
	
	Before giving the proof of Theorem \ref{thm:hhperf2}, we need the following:
	\begin{lem}
		Let $A$ be a curved algebra. The functor
		\begin{align*}
			(-)^\vee=\Homdg_{A}(-,A):\Tw(A)^\cop&\to \Tw(A^\op)
		\end{align*}
		is an equivalence of dg categories. 
	\end{lem}
	\begin{proof}
		Let $M\in \Tw(A)^\cop$, then by definition $M^\#$ is a finitely generated graded free $A^\#$-module. Every such $A$-module is mapped to a finitely generated twisted left $A$-module under the functor $(-)^\vee$. An inverse of $(-)^\vee$ is given by $\Homdg_{A^\op}(-,A): \Tw(A^\op)\to \Tw(A)^\cop$ thus $(-)^\vee$ must be an equivalence. 
	\end{proof}
	\begin{proof}[Proof of Theorem \ref{thm:hhperf2}]
		We first compute $\HHc(\Tw(A))$. Let
		\begin{align*}
			\Phi:\DGMod[A]{A}&\to \DGMod{\Tw(A^\op)\otimes \Tw(A)}\\
			M&\mapsto \Homdg_{A^e}(-,M)
		\end{align*}
		be the functor used in the proof of Theorem \ref{T:twbim} and define the dg equivalence
		\begin{align*}
			\Psi=\left((-)^\vee\right)^\cop\otimes \id_{\Tw(A)^\cop}:\Tw(A)\otimes \Tw(A)^\cop\to \Tw(A^\op)^\cop\otimes \Tw(A)^\cop.
		\end{align*}
		By precomposing with $\Psi$, we obtain the functor
		\begin{align*}
			\Psi^\ast\Phi:\DGMod[A]{A}\to \DGMod[\Tw(A)]{\Tw(A)},
		\end{align*}
		defined by sending the $A$-bimodule $M$ to the $\Tw(A)$-bimodule $\Psi^\ast\Phi(M)=[K\otimes L\mapsto \Homdg_{A^e}(K^\vee\otimes L,M)]$. There is a natural isomorphism
		\begin{equation*}
			\Homdg_{A^e}(K^\vee\otimes L,M)\cong \Homdg_{A}(L,\Homdg_{A^\op}(K^\vee,M))\cong \Homdg_{A}(L,K\otimes_A M),
		\end{equation*}
		where the first isomorphism is the $\otimes$-Hom adjunction and the second isomorphism follows from dualizability of finite generated twisted modules. 
	
	Under the isomorphism above, we identify $\Psi^\ast\Phi(M)$ with the functor $[K\otimes L\mapsto \Homdg_{A}(L,K\otimes_A M)]$. Clearly, $\Psi^\ast\Phi(A)$ is the diagonal module $\Tw(A)$, and since $\Psi^*$ and $\Phi$ preserve all weak equivalences, the derived functor sends $A\in \DII(A,A)$ to $\Tw(A)\in \D(\Tw(A)^\cop\otimes \Tw(A))$.
		
		We  then compute
		\begin{align*}
			\HHc(\PerfII(A))&\simeq \HHc(\Tw(A))\\
			&\simeq \RHom_{\Tw(A)^\op \otimes \Tw(A)}(\Tw(A), \Tw(A)) \\
			&\simeq \RHom^\text{II}_{A,A}(A,A) \\
			&\simeq \HHIIc(A).
		\end{align*}
		The first quasi-isomorphism follows from derived Morita equivalence, the third is a consequence of Theorem \ref{T:twbim} together with the computation above. 
		
		This proves the first quasi-isomorphism in the statement of the theorem. The second follows from Corollary \ref{comp:bimod} and \cite[Corollary 8.2]{toe06}.
	\end{proof}
			In particular if $A$ and $E$ are II-Morita equivalent in the sense that $\PerfII(A)$ and $\PerfII(E)$ are quasi-equivalent dg categories, then $\HHII(A) \simeq \HHII(E)$. (There is no need to assume the equivalence of categories is induced by a morphism or $(A, E)$-bimodule.)
			
			In particular, any MC equivalence induces a quasi-isomorphism on Hochschild cohomology of the second kind. 
						\begin{cor}\label{C:mchh}
				An MC equivalence $F \colon A \to B$ between curved algebras induces a quasi-isomorphism $\HHIIc(A) \simeq \HHIIc(B)$ of dg algebras. 
			\end{cor}
			\begin{proof}
				This is a direct consequence of Theorem \ref{thm:hhperf2} and Proposition \ref{prop:mcd2}.
			\end{proof}
	\begin{cor}
		For any curved algebra	$\HHIIc(A)$ is quasi-isomorphic to a $B_\infty$ algebra. In particular, 
	$\HHII(A)$ has the structure of a Gerstenhaber algebra.
	\end{cor}
	\begin{proof}
		This follows from the well-known $B_\infty$-struture on $\HHc$ of any dg category \cite{kel06a}.
	\end{proof}

\subsection{Koszul invariance}\label{sect:koszulinvariance}
We will now show that Koszul duality identifies suitable notions of Hochschild cohomology. 

By $\HHctr(C)$ for a  pc algebra $C$ we denote the Hochschild cohomology computed in the contraderived category of $C$-bimodules, i.e.\ the cohomology of the complex $\HHctr(C) = \RHom_{C \otimes C^\op}(C,C)$.
This sometimes called \emph{coHochschild cohomology} of the  dual coalgebra $C^*$. 
\begin{cor}\label{C:computeOmega}
	Let $C$  be pc algebra. Then $\HHctr(C) \simeq \HHIIc(\Omega C)$.
\end{cor}
\begin{proof}
	By Theorem \ref{T:monoidalCobar},
	there is an equivalence between
	$\DII(\Omega C^\op \otimes \Omega C)$ and $\DII(\Omega (C^\op \otimes C)) \cong \Dctr(C^\op \otimes C)$.
	Thus, the result follows if we can identify the diagonal bimodules.
	
	We consider  $\hat G C = (\Omega C \otimes C^* \otimes \Omega C)^{[\xi]}$; then the the natural comparison map $\hat G C \to \Omega C$ is identified with the morphism $W_{\Omega C}\to\Omega C$ where $W_{\Omega C}$ is the $\Omega C$-bimodule constructed in Lemma \ref{lem:curved}. Since the kernel of this surjective map is contractible by Lemma \ref{lem:curved}, it is a weak equivalence in $\DGModII[\Omega C]{\Omega C}$ by Proposition \ref{prop:useful}.			
\end{proof}
The first part of the next result  appears in \cite[Proposition 2.4]{kel21a}. 
\begin{cor}\label{C:computeB}
	For any dg algebra $A$, there are quasi-isomorphisms of algebras
	\[
	\HHc(A) \simeq \HHctr(\hB A)
	\]
	and for any curved algebra $A$
	\[
	\HHIIc(A) \simeq \HHctr(\cB A).
	\]
\end{cor}
\begin{proof}
	We prove the second statement. By Corollary \ref{C:computeOmega}, it suffices to show that $\HHIIc(A) \simeq \HHIIc(\Omega \check B A)$. But this follows from Theorem \ref{C:mchh} as $ \Omega \check B A \to A$ is an MC equivalence, see Example \ref{exmp:cobarbar}.
\end{proof}

For convenience we formulate the next corollary in terms of coalgebras.

\begin{cor}
	Let $C$ be a pseudocompact curved algebra. Then $\HHctr(C)$ is isomorphic to the 	Hochschild cohomology of its (dg enhanced) contraderived category.
\end{cor}
\begin{proof}
	We use the subscript $dg$ for dg enhancement again.
By Corollary \ref{C:computeOmega}, Theorem \ref{thm:hhperf2} and	Theorem \ref{T:KdCobar} in order we have
\begin{align*}
\HHctr(C) &\simeq \HHIIc(\Omega C) \simeq \HHc(\DIIdg(\Omega C)) \simeq \HHc(\Dctrdg(C)).  \qedhere
\end{align*}
\end{proof}
Equivalently, the coHochschild homology of an arbitrary curved coalgebra is quasi-isomorphic to the Hochschild cohomology of its coderived category.

	\section{Examples}\label{sect:examples}
	
	\subsection{Complex manifolds}
		Let $X$ be a complex manifold and $A(X):=(\mathcal{A}^{0,*}(X), \bar \partial)$ its Dolbeault-algebra. 
	It is well-known \cite{blo10, chl21} that the bounded derived category $\DCoh^b(X)$ of coherent sheaves on $X$ is equivalent to the derived category of sheaves on $X$ with bounded  coherent cohomology and, since $X$ is smooth, the latter coincides with the derived category of perfect complexes of sheaves on $X$. 
		
		We recall the following result relating $A(X)$ and $\Perf(X)$:
	\begin{thm}\label{thm:perfiidolbeault}
		Let $X$ be a complex manifold.
		$\PerfII(A(X))$ is quasi-equivalent to $\Perf(X)$.
	\end{thm}
	\begin{proof}
		This is essentially the main result of Block \cite{blo10}, who considers the category of cohesive modules over the Dolbeault algebra. For the equivalence of this category with $\PerfII(A(X))$ see \cite[Corollary 3.1]{Holstein5}, where $\PerfII$ is called the category of perfect twisted modules.
	\end{proof}
	
	\begin{thm}\label{thm:hhdolbeault}
		Let $X$ be a complex manifold. We have a natural isomorphism
		\[
		\HHII(A(X))\cong \HH(\Perf(X)).
		\]
		In particular we can compute $\HHII(A(X)) \cong H^*(X, \wedge^* \mathcal T_X)$.
	\end{thm}
	\begin{proof}
		This is just an application of Theorem \ref{thm:hhperf2}.
		As $A$ is graded commutative, $A^\op\cong A$ and $\Perf(A)\simeq (\PerfII(A))^\op$.
		
		The computation now follows from the standard result that $\HH(\DCoh^b(X))$ is computed by $H^*(X, \wedge^* \mathcal T_X)$ which follows by combining \cite[Corollary 4.2]{calduararu2005mukai} and \cite[Corollary 8.12]{toen2007homotopy}
	\end{proof}	
	
	\subsection{Topological spaces}
	In the topological setting we consider $\infty$-local systems on a topological space $X$. To be precise, we denote by $\LC(X)$ the dg category of fibrant cofibrant cohomologically locally constant (clc)  sheaves of complexes over $k$ whose cohomology sheaves have finite-dimensional fibers. 
	We recall the following:
	\begin{prop}\label{prop:perf2ls}
		Let $X$ be a connected locally contractible topological space and $C^*(X)$ its normalized singular cochain algebra with coefficients in $k$.
		Then $\PerfII(C^*(X)) \cong \LC(X)$.
		
		Let $k = \mathbb R$, $X$ be a connected manifold and $\mathcal A^*(X)$ be its de Rham algebra. Then $\PerfII(\mathcal A^*(X)) \cong \LC(X)$.
	\end{prop}
	\begin{proof}
		For the first claim, it follows from \cite[Theorem 8.4]{chl21} that $\LC(X)$ is given by finitely generated twisted modules over $C^*(X)$ since
		finitely generated twisted modules correspond to clc sheaves whose fibers are bounded and finite dimensional in each degree. As $\LC(X)$ is idempotent complete it follows that it is quasi-equivalent to $\PerfII(C^*(X))$.
		
		The second claim is \cite[Theorem 8.1]{chl21}.
	\end{proof}
	\begin{rem}\label{rmk:contra}
		We may go beyond $\infty$-local systems with perfect fibers  if we consider $C^*(X)$ as a pc algebra. 
		We consider $\Dctr(C^*(X))$ as the homotopy  category of modules over $C^*(X)$ whose underlying graded is of the form $C^*(X, V)$ for some chain complex $V$.
		Then there is an equivalence $
		\Dctr(C^*(X)) \cong \LC^\infty(X)$, where the right hand side is the category of potentially infinite-dimensional local systems \cite[Theorem 8.4]{chl21}.
		
		But $\Dctr(C^*(X))$ is different from $\DII(C^*(X))$, the natural map $\DII(C^*(X)) \to \Dctr(C^*(X))$ has in its image all  $\infty$-local systems that arise as colimits of $\infty$-local systems with perfect fibers. 
		For $X = S^1$ this is clearly not essentially surjective as the regular representation of $k[\Z]$ is not a colimit of finite-dimensional ones.
		
		Note that if $\Omega X$ is of finite type then we may use $\LC^\infty(X) \simeq D(C_*\Omega X)$ to see that $\LC^\infty(X)$ is compactly generated by finite-dimensional local systems and equivalent to $\DII(C^*(X))$.
		
		In particular we can deduce that $\HHII(C^*X) \simeq \HH(C_*\Omega X)$ in this case.
	\end{rem}

	\begin{cor}
		Let $X$ be a connected locally contractible topological space. Then $\HHII(C^*(X)) \cong \HH(\LC(X)))$.
		
		If  $X$ is moreover a manifold, then $\HHII(\mathcal A^*(X)) \cong \HH(\LC(X)))$.
	\end{cor}
	
	\begin{proof}
		This is a direct consequence of Theorem \ref{thm:hhperf2} and Proposition \ref{prop:perf2ls}.	
	\end{proof}

		\subsection{Matrix Factorizations}
	We now turn to a curved example. Let $R$ be a regular commutative algebra over a field $k$ of characteristic 0, and let $w \in R$. 
	Interpreting $R$ as a $\Z/2$ graded algebra concentrated in even degrees we may consider $w$ as curvature and define a curved ring $R_w = (R, 0, w)$.
	
	Then the dg category of matrix factorizations $\MF(R,w)$ may be defined as the idempotent completion of the category of curved modules over $R_w$ such that the underlying graded module is finitely generated and projective in each degree. 
	Its homotopy category is known to be equivalent to the derived category of singularities $\D^b_{coh}(Z)/\Perf(Z)$ where $Z = w^{-1}(0) \subset \operatorname{Spec}(R)$ is the zero locus with singular locus $\operatorname{crit}(w)$ \cite{Orlov03}.
	
	For a $\Z/2$-graded curved algebra we may define its $\Z/2$-graded categories of twisted modules and perfect complexes of the  second kind $\Perf^{\textrm{II}}_{\Z/2}$ by just changing the grading in our definitions. 
	
	With these definitions we obtain the  following lemma:
	\begin{lem}\label{lem:mf}
		There is a quasi-equivalence $\MF(R,w) \simeq \Perf^{\mathrm{II}}_{\Z/2}(R_w)$ of $\Z/2$-graded dg categories.
	\end{lem}
	\begin{proof}
		The definitions of $\MF(R, w)$ and $\Perf^{\mathrm{II}}_{\Z/2}(R_w)$ agree except that matrix factorizations are built out of projective modules for $R^\#$,  the underlying graded  of $R$, rather than free modules.
		But any curved $R_w$ module whose underlying graded module $P$ is finitely generated projective is a direct summand of a finitely generated twisted module, i.e.\ a curved module whose underlying graded is finitely generated free over $R^\#$.
		
		To show this, pick a finitely generated projective $R^\#$-module $L$ such that $P \oplus L$ is free and consider $G(L)$ the free curved $R$-module on $L$ which has elements formal sums $\ell + d\ell$ and differential $d(\ell + d\ell) = h\ell + d\ell$, cf.\ \cite[Section 3.6]{pos11}.
		Then the underlying graded of $G(L)$ is $L\oplus L[-1]$ and $P \oplus P[-1] \oplus G(L)$ is the desired $R_w$-module, proving the assertion.
	\end{proof}
	
	We define Hochschild cohomology of the  second kind of $R_w$ as in the $\Z$-graded case, noting that 
	$\HHII(R_w) = \RHomII_{R_w \otimes R_w^{\op}} (R_w, R_w)$ is now a $\Z/2$-graded complex.
	
	\begin{prop}\label{prop:MF}
		With $(R,w)$ as above we have $\HHII(R_w) \cong \HH(\MF(R,w))$.
		
		If $\mathrm{crit}(w) \subset w^{-1}(0)$ then $\HHII(R_w)$ can be computed explicitly as the cohomology of the complex $(\wedge^* \operatorname{Der}(R)[1], \iota_{dw})$ where the differential is given by contraction with $dw$.
	\end{prop}
	Note that the first statement does not make any assumptions on $(R,w)$, neither isolated singularities nor even the assumption $\mathrm{crit}(w) \subset w^{-1}(0)$.
	\begin{proof}
		Theorem \ref{thm:hhperf2} and Lemma \ref{lem:mf} immediately give \[\HHII(R_w) \cong \HH(\MF(R,w)).\]
		The explicit formula for $\HHII(R_w)$ comes from \cite[Theorem 3.1]{lin2013global} who write $\iota_{df}$ equivalently as the Schouten bracket $[f,-]$.
	\end{proof}
	
	\begin{rem}\label{rem:mfhhppagree}
		Note that $\HH(\MF(R,w))$ has also been computed in terms of compactly supported Hochschild cohomology (cf.\ Remark \ref{remark:comparison}). 
		For isolated singularities of $w^{-1}(0)$ this is \cite{caldararu2013curved, tu2014matrix} 
		and if $\mathrm{crit}(w) \subset w^{-1}(0)$ this is
		\cite[Section 4.8]{pp12}.
		They show $\HH(\MF(R,w)) = 
		\HH_{\mathrm{PP}}(R_w)$. 
		This is further used for example in \cite{efimov2018cyclic}.
		
		It is notable that in this important case our definitions agree with the older definition that is different in many cases. We consider the reason for this in Section \ref{sect:hhpp}.
		
		If there are critical values other than $0$ then  $$\HH_{\mathrm{PP}}(R_w) \cong \oplus_{c \in C} \HH(\MF(R, w-c))$$
		where the sum is over the set $C$ of critical values, see \cite[Section 4.10]{pp12}.
		Together with the Proposition, \ref{prop:MF}, this gives
		$$\HH_{\mathrm{PP}}(R_w) \cong \oplus_{c \in C} \HHII(R_{ w-c}) \cong (\wedge^* \operatorname{Der}(R)[1], \iota_{dw}).$$	
	\end{rem}
	

	\section{Comparison results}
	
	\subsection{Comparison of model structures for bimodules}
After noting the different model structures on $(A,E)$-bimodules and on $A\otimes E^\op$-modules, we now observe that in many cases of interest the constructions do agree.

	Consider 
	the natural functor 
	\[
	i \colon \Perf(\PerfII(A)\otimes \PerfII(E))\to  \PerfII(A\otimes E)
	\]
	sending a pair $M, N$ of perfect modules of the second kind to $M \otimes N$.
	
	We note that this is not in general an equivalence and consider the following condition. 
	
\begin{con}\label{con:star}
We say that a pair of curved algebras $(A, E)$ satisfies Condition \ref{con:star} if 	\[
	i \colon \Perf(\PerfII(A^\op)\otimes \PerfII(E))\to  \PerfII(A^\op\otimes E)	
	\]
	is a quasi-equivalence. We say $A$ satisfies Condition \ref{con:star} if $(A, A)$ does.
	\end{con}	
	
	\begin{rem}
		A similar condition appears in the work of Polishchuk-Positselski when investigating Morita invariance of $\HH_{\mathrm{PP}}$, see e.g.\ \cite[Theorem 3.5.C]{pp12}.
		The functor $i$ (and its failure to be an equivalence) is also considered in \cite[Section 2]{briggs2024hochschildhomologycurvedalgebras}.
	\end{rem}

	\begin{rem}
		The functor $i$ is not always a quasi-equivalence as illustrated by Example \ref{eg:counterexample to condition 1}.
		The dg category $\PerfII(A\otimes E)$ should be viewed as a kind of \emph{completed} tensor product of the categories $\PerfII(A)$ and $\PerfII(E)$. 
		There are, however, important situations when this completion is extraneous.
	\end{rem}
	\begin{exmp}\label{eg:counterexample to condition 1}
		The functor $i$ associated with a pair $(A,E)$ is generally not essentially surjective:
		\begin{itemize}
			\item[(1)] Let $A$ be the $\mathbb{Z}/2\mathbb{Z}$-graded curved algebra $(k,0,1)$. A finitely generated twisted $A$-module $M$ corresponds to a pair of $k$-linear maps 
			\begin{align*}
				f:M\rightleftarrows N:g
			\end{align*}
			such that $fg=\operatorname{id}$ and $gf=\operatorname{id}$. Every such $A$-module is contractible, thus $\PerfII(A)\simeq 0$. The enveloping algebra $A^e$ is the $\mathbb{Z}/2\mathbb{Z}$-graded (curved) algebra $(k,0,0)$. The category $\PerfII(A^e)$ is equivalent to the category of finite-dimensional super vector spaces and is therefore not Morita equivalent to $\PerfII(A^\op)\otimes \PerfII(A)$
			\item[(2)] Let $A=\Lambda(\epsilon)$ with $|\epsilon|=-1$ and $E=k[x,x^{-1}]$ with $|x|=2$. For dg algebras concentrated in non-positive degrees, the extended bar construction coincides with the ordinary one, which by Koszul duality implies that the categories $\PerfII(A)$ and $\Perf(A)$ are equivalent.
			As $E$ is a graded field, $K(E)\simeq \DII(E)\simeq \D(E)$ thus $\PerfII(E)\simeq \Perf(E)$.
			At the same time, it can easily be checked that $\PerfII(A^\op\otimes E)\not\simeq \Perf(A^\op\otimes E)$. For instance, the $A^\op\otimes E$-module $(A\otimes E)^{[\epsilon\otimes x]}\in \PerfII(A^\op\otimes E)$ is acyclic but non-contractible and so it represents the zero object in $\Perf(A^\op\otimes E)$ but a non-zero object in $\PerfII(A^\op\otimes E)$.
		\end{itemize}
	\end{exmp}
	We have the following consequence of Condition \ref{con:star}.
	\begin{thm}\label{thm:condition1compare}
		If a pair of curved algebras $A$ and $E$ satisfies Condition \ref{con:star}, the adjunction
		$$\operatorname{id}:\DGMod[A]{E} \rightleftarrows \DGMod{A^\op \otimes E}:\operatorname{id}$$
		is a Quillen equivalence  
	\end{thm}
	\begin{proof}
		The Quillen adjunction was noted in Remark \ref{rem:naivebarquillen}. 
		By construction $\DII(A^\op \otimes E)$ is compactly generated by $\PerfII(A^\op \otimes E)$, which by assumption is equivalent to the idempotent closure of $\PerfII(A)^\op \otimes \PerfII(E)$.
		But the latter category compactly generates $\DII(A^\op \otimes E)$, proving the adjunction is an equivalence on homotopy categories.
	\end{proof}

	We conclude the section with two results that are not needed in the rest of the paper, but the first is a useful tool for working with curved algebras and the second clarifies the behaviour of $\PerfII$.
	
	We first note the following lemma:
	\begin{lem}	\label{L:perf4}
		For any curved algebra $A$ the Yoneda embedding 
		$\PerfII(\cat{A}) \to \PerfII(\PerfII(\cat{A}))$ is a quasi-equivalence.
	\end{lem}
	
	\begin{proof}
		By definition $\PerfII(\PerfII(A))$ is the homotopy idempotent closure of $\Tw(\PerfII(A))$. But as $\PerfII(A)$ is closed under finite direct sums and MC twists it is pre-triangulated in the sense of \cite{bk91}. 
		Thus  $\Tw(\PerfII(A))$ is equivalent as a dg category to $\PerfII(A)$ by \cite[Proposition 3.1]{bk91}.
		Hence it is already idempotent complete and $\PerfII(\PerfII(A)) \cong \Tw(\PerfII(A)) \cong \PerfII(A)$. 
	\end{proof}

	\begin{lem}\label{lem:uncurved}
	Given a curved algebra $A$, there is a dg algebra $A'$ such that $\PerfII(A)$ and $\PerfII(A')$ are equivalent dg categories.
\end{lem}
We may say $A'$ is II-Morita equivalent to $A$.
\begin{proof}	
	Let $w\in A^2$ be the curvature element. Let $V:=k\oplus k[1]$ and consider the $A$-module $M:=A\oplus A[1]\cong A\otimes V$ with the differential $d_M$ given on $1\otimes V\subset M$ by the $2\times 2$ matrix $x_w=-\begin{pmatrix}0&1\\w&0\end{pmatrix}$. 
	We have $d_M^2=\begin{pmatrix}w&0\\0&w\end{pmatrix}$ and so $M$ is a curved (perfect) $A$-module. It follows that $x_w$ is an MC element in the curved algebra $A\otimes \End(V)$. Set $A':=\End_A(M)\cong A\otimes \End(V)^{x_w}$, the twist of $A\otimes \End(V)$ by $x_w$.
	We clearly have an equivalence $\PerfII(A)\to\PerfII(A\otimes \End(V))$ induced by the usual prescription $?\mapsto \Hom_A(?,M)$. 
	
	Since $A'$ is isomorphic to $A\otimes \End(V)$ as a curved algebra, it follows that $\PerfII(A)$ and $\PerfII(A')$ are equivalent as dg categories. 
\end{proof}

\begin{rem}
	It is easy to see that the dg algebra $A'$ constructed above, is acyclic (and so, the twisted module $M$ is homotopically trivial). Indeed, the identity element $\begin{pmatrix}1&0\\0&1\end{pmatrix}$ is the coboundary of the element $\begin{pmatrix}0&0\\-1&0\end{pmatrix}$ in $A'$. Note that the ordinary derived category of $A'$ is, of course, trivial. 
\end{rem}

The following statement for $\PerfII$ is analogous to the corresponding statement for $\Perf$.

	\begin{prop}\label{L:perfperf}
		For curved algebras $A, B$ there is a quasi-equivalence
		\[
		\PerfII(A \otimes B) \simeq \PerfII(\PerfII(A) \otimes \PerfII(B))
		\]
	\end{prop}
	\begin{proof}
		We consider the case that $A$ and $B$ are dg algebras first.
		Consider the following diagram
		\[
		\begin{tikzcd}[sep=large]
			\PerfII(A \otimes B) & \PerfII(\PerfII(A) \otimes \PerfII(A)^{\op}) \\
			& \PerfII(\PerfII(A \otimes B))
			\arrow["{(h^A\otimes h^B)_!}", from=1-1, to=1-2]
			\arrow["{i_!}", from=1-2, to=2-2]
			\arrow["{h^{A\otimes B}_!}", from=1-1, to=2-2]
		\end{tikzcd}
		\]
		Here $h^{A \otimes B}_!$ is a quasi-equivalence by Lemma \ref{L:perf4}.
		as the Yoneda embedding induces an equivalence on $\PerfII$ by Lemma \ref{L:perf4}.
		The functor $i: \PerfII(A)\otimes \PerfII(B) \to \PerfII(A \otimes B)$ is fully faithful (not just quasi-fully faithful)
		and thus $i_!$ is also fully faithful.
		Moreover $i_!$ is quasi-essentially surjective as it has the same image as the quasi-equivalence $h_!^{A\otimes B}$.
		It follows that $(h^A \otimes h^B)_!$  is a quasi-equivalence by 2-out-of-3. 
	
	Assume now that $A, B$ are curved algebras.
		We have an equivalence of categories $\PerfII(A\otimes B)\cong \PerfII(A'\otimes B')$ (where $A'$ and $B'$ are dg algebras II-Morita equivalent to $A$ and $B$ as constructed in Lemma \ref{lem:uncurved}) since $A'\otimes B'$ is isomorphic to the tensor product of $A\otimes B$ and a $4\times 4$ matrix algebra. 
		Similarly $\PerfII(\PerfII(A')\otimes \PerfII(B'))$ is quasi-equivalent to $\PerfII(\PerfII(A)\otimes \PerfII(B))$ which reduces the question to the uncurved case.
	\end{proof}

	\subsection{Examples}
	
	One reason to include this discussion of the relation between the model categories $\DGMod[A]{E}$ and $\DGMod{A^\op \otimes B}$ is that the three classes of examples from Section \ref{sect:examples} satisfy Condition \ref{con:star} under some additional assumptions.
	
	We first observe that Condition \ref{con:star} is satisfied whenever $A$ is cofibrant in $\cuAlg$ by Corollary \ref{C:bimoduleovercofibrant}
	, but this is very much not the case for the examples of geometric and topological interest that we are considering now.
	
	Let $X$ be a complex projective manifold and with bounded derived category $\DCoh^b(X)$ of coherent sheaves.
	We will consider its dg model $\Perf(X)$ formed by taking Dolbeault resolutions of coherent sheaves. 
	According to \cite[Corollary 3.1.8]{bvdb03}, the latter is equivalent to the dg category of dg modules over some dg algebra (the endomorphism algebra of a generator of $\Perf(X)$).

	\begin{prop}\label{prop:tensor}
		Let $X$ be a smooth projective variety. 
		The categories $\Perf(X\times X)$ and $\Perf(X)\otimes\Perf(X)$ are Morita equivalent.
	\end{prop}
	\begin{proof}
		The statement is well-known and follows from 
		a very general result \cite[Theorem 1.2(1)]{bfn10} 
		valid for perfect derived stacks, not merely for complex projective manifolds. The argument goes back to \cite{bvdb03}.
		We sketch a proof for the reader's convenience. 
		
		Let $\mathcal{F}$ be a complex of sheaves representing a generator of $\Perf(X)$ and $B=R\End(\mathcal{F})$ be its endomorphism dg algebra. 
		Then the external tensor product $\mathcal{F}\boxtimes \mathcal{F}$ is a generator of $\Perf(X\times X)$ and $B\otimes B\simeq R\End(\mathcal{F}\boxtimes \mathcal{F})$. So, $\Perf(X)$ is quasi-equivalent to $\Perf(B)$ and $\Perf(X\times X)$ is quasi-equivalent to $B\otimes B$. 
		Thus we have reduced the question to that of a perfect derived category of a dg algebra. Since $\Perf(B\otimes B)$ is Morita equivalent to $\Perf(B)\otimes\Perf(B)$ (for any dg algebra $B$), the desired claim follows. 
	\end{proof}	

We can deduce that Condition \ref{con:star} holds in the algebraic case (but not for non-algebraic smooth manifolds).
	
	\begin{thm}\label{thm:hhdolbeault}
		Let $X$ be a smooth projective variety. 
		The Dolbeault algebra $A(X)$ satisfies Condition \ref{con:star}.
	\end{thm}
	\begin{proof}
		As $A$ is graded commutative $A^\op\cong A$ and $\Perf(A)\simeq (\PerfII(A))^\cop$.
		
		From Proposition \ref{prop:tensor} 
		together with Theorem \ref{thm:perfiidolbeault} we obtain
		that 
		$\PerfII(A(X)) \otimes \PerfII(A(X))$ is Morita equivalent $\PerfII(A(X \times X))$.	
		It remains to show that $\PerfII(A(X \times X))$ is Morita equivalent to $\PerfII(A(X) \otimes A(X))$.
	It follows from \cite[Theorem 51.6]{Treves67} that there is an isomorphism of dg algebras between $A(X\times X)$ and the completed tensor product $A(X)\hat \otimes A(X)$ and thus we have Morita equivalences
$$\PerfII(A(X))\otimes \PerfII (A(X)) \simeq 
\PerfII(A(X\times X)) \simeq \PerfII(A(X) \hat \otimes A(X)).$$
	The composition of the equivalences naturally factors as
	$$\PerfII(A(X)) \otimes \PerfII(A(X)) \to \PerfII(A(X) \otimes A(X))) \to \PerfII(A(X) \hat \otimes A(X)) $$
	which shows that the second map induces a quasi-essentially surjective map on derived categories.
	As the induced map is also quasi-fully faithful by the Künneth theorem it must be a quasi-equivalence and putting all our equivalences together we obtain the desired Morita equivalence
	$\PerfII(A(X) \otimes A(X)) \simeq \PerfII(A(X)) \otimes \PerfII(A(X))$.
	\end{proof}

Similar arguments apply in the topological setting.
	\begin{prop}\label{prop:lcmult}
		Let $X$ be a CW complex with finitely many cells in each degree. Then $\LC(X 
		\times X)$ and $\LC(X) \otimes \LC(X)$ are Morita equivalent.
	\end{prop}
	\begin{proof}
		We first show the natural functor given by the exterior product is essentially surjective. Let $G$ be the fundamental group of $X$. As $LC(X)$ is classically generated by irreducible local systems, it suffices to show that any finite-dimensional irreducible module over $k[G] \otimes k[G]$ is a summand of a tensor product of irreducible $k[G]$-modules. This follows from Lemma \ref{lem:separablealgebras}.
		
		It remains to be shown the functor is quasi-fully faithful. It suffices to check on generators, i.e.\ on local systems. 
		We thus have to compare
		$\RHom(L, L') \otimes \RHom(M, M')$ with $\RHom(L \boxtimes L, M \boxtimes M')$ for local systems $L,L', M, M'$ on $X$.
		But this is equivalent to showing
		$$C^*(X, L' \otimes L^*) \otimes C^*(X, M' \otimes M^*) \simeq C^*(X \times X, (L' \otimes L^*)\boxtimes (M' \otimes M^*)).$$
		Thus the result follows from the K\"unneth theorem with local coefficients \cite[Theorem 1.7]{Greenblatt06}. Here we use the finiteness assumption. 
	\end{proof}
	\begin{cor}
		Let $X$ be a connected locally contractible topological space that has the homotopy type of a CW complex with finitely many cells in each degree. Then $C^*(X)$ satisfies Condition \ref{con:star}.
	\end{cor}

	\begin{proof}
		By Proposition \ref{prop:perf2ls} we have $\PerfII(C^*(X)) \simeq \LC(X)$, so the result follows if we can verify that
		$\PerfII(C^*(X))\otimes \PerfII(C^*(X)) $ and $\PerfII(C^*(X) \otimes C^*(X)^{\op})$ are Morita equivalent.		
		From Proposition \ref{prop:lcmult} we have $\PerfII(C^*(X)) \otimes \PerfII(C^*(X))$ is Morita equivalent to $\PerfII(C^*(X\times X))$.
		
Thus it remains to show that there is a Morita equvialence between $\PerfII( C^*(X \times X)) $ and $\PerfII(C^*(X)  \otimes C^*(X))$ to complete the proof.
		
		We know $D(\Omega(C_*(X \times X)) \simeq \LC^\infty(X \times X)$. 
		 As $\Omega$ is quasi-strong monoidal \cite{holstein2025enriched} we also have $D(\Omega (C_*X \otimes C_*(X))) \simeq D(\Omega(C_*X) \otimes \Omega(C_*(X))) \simeq \LC^\infty(X\times X)$.
		Thus by Koszul duality we have an equivalence
		$\Dctr(C^*(X \times X)) \simeq \Dctr(C^*(X) \otimes C^*(X))$, cf.\ Remark \ref{rmk:contra}. 
		It remains to identify the subcategories of finitely generated twisted modules.
		Replacing singular chains by reduced singular chains we may assume $C^0(X) = k$.
		The equivalence $\Dctr(C^*(X \times X)) \simeq \Dctr(C^*(X) \otimes C^*(X))$ is compatible with the maps $M \mapsto M \otimes_A k$ for $A = C^*(X\times X)$ respectively $C^*(X)\otimes C^*(X)$. 
		
		Thus it identifies the modules $M$ such that $M \otimes_A k$ is finite-dimensional. But these modules are exactly the finitely generated twisted modules.
		One inclusion is clear and given $M \otimes_A k$ finite-dimensional we can use \cite[Proposition 6.4]{Holstein5} to show $M$ is equivalent to a twisted finitely generated module.		
	\end{proof}

We finally consider the case of matrix factorization. Here the desired result is an immediate consequence of the Thom-Sebastiani-Theorem.

		Let $k$ have characteristic 0 and let $R$ be a regular commutative $k$-algebra with $w \in R$.
		
	\begin{prop}
	Let $(R,w)$ be as above  such that $\operatorname{crit}(w) \subset w^{-1}(0)$. Then the curved $\Z/2$-graded algebra $R_w$  satisfies Condition \ref{con:star}.
	\end{prop}
	\begin{proof}
		The Thom-Sebastiani theorem \cite[Theorem 4.1.3]{Preygel11} says 
		$\MF(R, w) \otimes \MF(R, w) \cong \MF(R \otimes R, w \otimes 1 + 1\otimes w)$. In terms of curved algebras, this is exactly saying $\MF(R_w) \otimes \MF(R_w) = \MF(R_w \otimes R_w)$. 
	\end{proof}
	
	\subsection{Comparison with Polishchuk-Positselski's $\HH_{\mathrm{PP}}$}\label{sect:hhpp}
	
	Motivated by the agreement of Hochschild cohomology for matrix factorizations observed in Remark \ref{rem:mfhhppagree},
	we give some natural sufficient conditions for an isomorphism $\HH_{\mathrm{PP}}(A) \cong \HHII(A)$.
	\begin{thm}\label{thm:comp HH_PP and HH2}
		Let a curved algebra $A$ satisfy the following three conditions:
		\begin{enumerate}
			\item $(A^\#)^e$ has finite homological dimension
			\item $(A^\#)^e$ is a coherent ring
			\item $A$ satisfies Condition \ref{con:star}. \label{con:three}
		\end{enumerate}
		Then there is a natural isomorphism $\HH_{\mathrm{PP}}(A) \cong \HHII(A)$. 
	\end{thm}
	Before presenting a proof, the following proposition is needed.
	\begin{prop}\label{comp:DII and Dco}
		Let $A$ be a curved algebra such that $A^\#$ is coherent and has finite homological dimension, then
		\begin{align}
			\DII(A)\simeq \Dco(A)\simeq \Dctr(A).
		\end{align} 
	\end{prop}
	\begin{proof}
		The second equivalence is proven in \cite[Theorem 3.6]{pos11} so it remains to show the first.
		
		 It follows from \cite[Corollary 8.20]{PositselskiStovicek2024}, that a set of compact generators $\textbf{C}\subset \Dco(A)$ can be obtained by taking all $A$-modules $M\in \Dco(A)$ such that $M^\#$ is a finitely presented $A^\#$-module.
		 
		 However, since $A^\#$ is coherent and has finite homological dimension, every such $A$-module $M$ admits a resolution of the form
		 $$0\to P_n\to ...\to P_0\to M,$$
		 where the rows are exact and every $P_i$ is a finitely generated graded projective $A$-module. See \cite[Theorem 3.6]{pos11} for a procedure on how to obtain such resolution. Clearly, $\operatorname{Tot}(P_n\to ...\to P_0)\in \PerfII(A)$ which implies that $\textbf{C}\subset \PerfII(A)$ thus the inclusion $\DII(A)\subset \Dco(A)$ must be an equivalence of categories.
	\end{proof}
	\begin{proof}[Proof of Theorem \ref{thm:comp HH_PP and HH2}]
   		It was proved in Proposistion \ref{comp:DII and Dco}, that under the assumption in Theorem \ref{thm:comp HH_PP and HH2}, there is a quasi-isomorphism
		$$\RHom_{\DII(A^e)}(A,A)\simeq \RHom_{\Dco(A)}(A,A).$$
		Moreover, since $(A^e)^\#$ has finite homological dimension, by \cite[Section 3.3]{pp12}
		$$\HH_{\mathrm{PP}}(A)\simeq \RHom_{\Dco(A)}(A,A).$$
		Finally Condition \ref{con:star} ensures we can apply Theorem \ref{thm:condition1compare} to deduce $\DGModII{A^e} \simeq \DGModII[A]{A}$ and putting everything together, we obtain
		\begin{align*}
			\HHIIc(A) &\simeq R\Hom_{\DGModII[A]{A}}(A,A)\\
			&\simeq R\Hom_{\DGModII{A^e}}(A,A)\\ &\simeq R\Hom_{\Dctr(A^e)}(A,A) \\
			&\simeq (\mathsf{B}'(A)\hat{\otimes}A)^{\xi}
		\end{align*}
		where the cohomology of the last term is $\HH_{\mathrm{PP}}(A)$.
	\end{proof}
Thus whenever the theorem applies we may compute $\HHII$ using $(\mathsf{B}'(A)\hat{\otimes}A)^{\xi}$ which is just the usual Hochschild complex except the direct product completion is replaced by a direct sum completion, making this complex much easier to work with than the extended bar construction.

The theorem also gives a conceptual explanation of what we observed in Remark \ref{rem:mfhhppagree}: The 2-periodization of a commutative regular ring $R$ satisfies the first two conditions, and if $w \in R$ satisfies $\operatorname{crit}(w) \subset w^{-1}(0)$ then the Thom-Sebastiani theorem ensures the final condition.
\begin{rem}
	Given a curved algebra $A$ the functor $i$ induces a fully faithful inclusion $\iota:\DII(A, A)\to \DII(A^e)$ of categories, even if $A$ does not satisfy Condition \ref{con:star}.
	
	Theorem \ref{thm:comp HH_PP and HH2} holds true even after relaxing condition \ref{con:three} to $A\in \DII(A^e)$ being in the essential image of $\iota$.
\end{rem}
	\bibliography{./hochschild_second_kind}

\setlength{\parindent}{0pt}
\end{document}